\documentclass[letterpaper,12pt]{article}

\usepackage[margin=1in,letterpaper]{geometry}

\usepackage[greek,english]{babel}
\usepackage[utf8]{inputenc}
\usepackage{amsmath,amssymb}
\usepackage{scalerel}
\usepackage{parskip}
\usepackage{graphicx}
\usepackage{mathrsfs}  
\usepackage{float} 
\usepackage{mathtools}
\usepackage{cancel}
\usepackage{tcolorbox}
\usepackage{amsmath, bm}
\usepackage{array}
\usepackage{amsthm}
\usepackage{csquotes}
\usepackage{hyperref}
\usepackage{amsthm}
\usepackage{quiver}
\usepackage{amsmath,amssymb}
\usepackage{scalerel}
\usepackage{parskip}
\usepackage{graphicx}
\usepackage{subcaption} 
\usepackage{mathrsfs}  
\usepackage{mathtools}
\usepackage{cancel}
\usepackage{tcolorbox}
\usepackage{amsmath, bm}
\usepackage{array}
\usepackage{amsthm}
\usepackage{csquotes}
\usepackage{hyperref}
\usepackage{amsthm}
\usepackage{quiver}

\DeclareUnicodeCharacter{0327}{}

\usepackage[style=ieee,backend=biber]{biblatex}

\DeclarePairedDelimiter\ceil{\lceil}{\rceil}
\DeclarePairedDelimiter\floor{\lfloor}{\rfloor}

\usepackage{tikz-cd}
\usetikzlibrary{cd}

\usepackage[shortlabels]{enumitem}

\usepackage{xcolor}

\makeatletter
\def\mathcolor#1#{\@mathcolor{#1}}
\def\@mathcolor#1#2#3{%
  \protect\leavevmode
  \begingroup
    \color#1{#2}#3%
  \endgroup
}
\makeatother

\theoremstyle{definition}
\newtheorem{exmp}{Example}[section]

\newcommand\biopencrossl{%
  \mathrel{\scalerel*{>\kern-.4\LMpt\joinrel\blacktriangleleft}{x}}}
\newcommand\biopencrossr{%
  \mathrel{\scalerel*{\blacktriangleright\joinrel\kern-.4\LMpt<}{x}}}

\newtheorem{theorem}{Theorem}

\newtheorem{corollary}{Corollary}

\newtheorem{lemma}[theorem]{Lemma}

\newtheorem{definition}{Definition}

\newtheorem{remark}{Remark}
\newtheorem{proposition}{Proposition}

\usepackage{mathrsfs}

\usepackage{xparse}

\NewDocumentCommand{\tens}{e{_^}}{%
  \mathbin{\mathop{\otimes}\displaylimits
    \IfValueT{#1}{_{#1}}
    \IfValueT{#2}{^{#2}}
  }%
}

\title{Assembly Addition Chains}
\author{Leroy Cronin* \\ Juan Carlos Morales Parra$\dagger$ \\ Keith Y. Patarroyo$\ddagger$  \\  
\small{University of Glasgow, Glasgow G11 6EW, United Kingdom} \\ 
\small{\url{lee.cronin@glasgow.ac.uk}}$*$ \\ \small{\url{juan.moralesparra@glasgow.ac.uk}}$\dagger$  \\ \small{\url{keith.patarroyo@glasgow.ac.uk}}$\ddagger$ 
\date{\today}}
\definecolor{lightblue}{RGB}{120,180,255}

\begin{document}

\maketitle
\begin{abstract} 
\noindent
In this paper we extend the notion of Addition Chains over $\mathbb{Z}^+$ to a general set $S$. We explain how the algebraic structure of Assembly Multi-Magma over the pairs $(S,BB\subset S)$ allows to define the concept of Addition Chain over $S$, called \textit{Assembly Addition Chains} of $S$ with Building Blocks $BB$. Analogously to the $\mathbb{Z}^+$ case, we introduce the concept of \textit{Optimal} Assembly Addition Chains over $S$ and prove lower and upper bounds for their lengths, similar to the bounds found by Schonhage for the $\mathbb{Z}^+$ case. In the general case the unit $1 \in \mathbb{Z}^+$ is replaced by the subset $BB$ and the mentioned bounds for the length of an Optimal Assembly Addition Chain of $ \mathscr{O} \in S$ are defined in terms of the \textit{size} of $\mathscr{O}$ (i.e. the number of Building Blocks required to construct $\mathscr{O}$). The main examples of $S$ that we consider through this papers are (i) $j$-Strings (Strings with an alphabeth of $j$ letters), (ii) Colored Connected Graphs and (iii) Colored Polyominoes. 
\end{abstract}

\section{Introduction and Motivation}

Addition Chains are \textit{representations} of positive numbers, constructed using 1 and the addition operation in the semi-group $(\mathbb{Z}^+,+)$. Formally, an Addition Chain (AC) for $n \in \mathbb{Z}^+$ is a finite increasing sequence of positive integers 
\begin{equation*}
    (a_0,a_1,\dots,a_r)
\end{equation*}
such that $a_0=1$, $a_r=n$ and for each $i \in \{1,2,\dots,r\}$ there exist $j,k \in \{0,1,\dots,i-1\}$ such that
\begin{equation}\label{adi}
    a_i=a_j+a_k.
\end{equation}
If we denote by $AC(n)$ the set of \textit{all}
the ACs of $n$, then each element of $AC(n)$ corresponds to a representation of $n$ as a recursive sum of positive integers starting from $1$. For example, for $n=7$ we have 
\begin{equation}\label{ac7}
\begin{split}
    AC(7)=\{ & (1,2,3,4,5,6,7),(1,2,3,4,5,7),(1,2,3,4,6,7),(1,2,3,4,7), (1,2,3,5,6,7), \\
    &(1,2,3,6,7), (1,2,3,5,7), (1,2,4,5,6,7),(1,2,4,5,7),(1,2,4,6,7)\},
\end{split}
\end{equation}
such that each of these finite sequences represents how the number $7$ could be constructed starting from $1$ and using recursively the binary operation of addition. \\\\
Among the set $AC(n)$, for each $n \in \mathbb{Z}^+$, there exists a particular subset of sequences exhibiting the shortest length. These are the so-called \textit{Optimal Addition Chains} (OAC), which we denote here by $OAC(n)$ for any $n \in \mathbb{Z}^+$, and they correspond consequently to the \textit{shortest representations} of a positive integer using $1$ and the recursive use of addition. For example, for $n=7$ then we have from (\ref{ac7}) that
\begin{equation}
  OAC(7)=\{(1,2,3,4,7),(1,2,3,5,7),(1,2,3,6,7),(1,2,4,5,7),(1,2,4,6,7)\}.  
\end{equation}
Given any $n \in \mathbb{Z}^{+}$, the length (i.e. the number of elements in the sequence) of any element in $OAC(n)$ is denoted by $\ell(n)$, e.g. $\ell(7)=5$.\\\\ 
The \textit{Optimal Addition Chain length function} 
\begin{equation}
\begin{split}
    \ell: \mathbb{Z}^+ & \to \mathbb{Z}^+ \\
    n & \mapsto \ell(n),
    \end{split}
\end{equation} 
besides the intrinsic natural mathematical interest, plays an important role in Computer Science and Complexity Theory. For example AC are used for \textit{Addition-Chain exponentiation}, which allows the computation of exponentiation with integer exponents by using a number of multiplications equal to the length of an AC for the exponent. Hence, in particular OAC are used to optimize the representation of any number in a computer (e.g. their binary representation). On the other hand $\ell(n)$ quantifies the \textit{complexity} of the integer $n$, in the sense that it corresponds to the minimal number of binary steps (addition) required to construct the number $n$ starting from the basic building block $1$. $\ell(n)$ could be understood as an additive restriction of the so-called Integer Complexity of an integer (see e.g. \cite{altman2015integer}). 
\\\\
We refer the interested reader to revisit the seminal work on AC by Scholz \cite{Scholz} and subsequent works done by Brauer \cite{brau}, Erd\"{o}s \cite{erdos1960remarks}, Stolarsky \cite{Stolarsky_1969}, Schonhage \cite{SCHONHAGE19751}, Clift \cite{clift}, Clift and Thruber \cite{THURBER2021112200}, among others. One of the best-known results was proved by Schonhage \cite{SCHONHAGE19751} and it shows that the length $\ell(n)$ for any $n \in \mathbb{Z}^+$ is bounded by       
\begin{equation}\label{ACineq}
     \log_2(n) + \log_2(H(n))-2.13 \leq \ell(n) \leq \log_2(n) + H(n)-1,
\end{equation}
where $H(n)$ is the Hamming weight of $n$ (i.e. the number of 1's in its binary representation). \\\\
The concept of AC requires only two ingredients: A set of Building Blocks (the $1$) and a binary gluing operation (the addition over $\mathbb{Z}^+$). So, in principle it is natural to extend the construction of AC (and consequently of OAC) over any set with a well-defined (i) set of Building Blocks and (ii) a binary gluing operation. \\\\
The goal of this paper is two-fold: First, we extend the concept of Addition Chains over $\mathbb{Z}^+$ to \textit{Assembly} Addition Chains over an arbitrary set $S$ (see Definition \ref{AAC}), showing how this generalization allows to define ACs for sets of (directed and undirected) strings, (colored connected) graphs and (colored) polyominoes, where the concepts of Building Blocks and binary gluing operation could be naturally described. Secondly, we prove for \textit{Optimal} Assembly Addition Chains an analogous result to (\ref{ACineq}), finding lower and upper bounds in terms of the \textit{size} of the object (i.e. the number of Building Block used to construct the object).\\\\
This extension of AC to the sets of strings, graphs, polyominoes, etc, is of profound relevance since, analogous to the AC for the elements of $\mathbb{Z}^+$, the Assembly Addition Chains provide \textit{representations} of the objects in terms of their recursive \textit{accessibility} from a fixed set of Building Blocks. In particular, the length of the shortest Assembly Addition Chain of an object (Assembly Index), i.e. one of its Optimal Assembly Addition Chains, quantifies the \textit{intrinsic complexity} of the object and it is intrinsically encoded in the object. \\\\ 
The notion of Multi-Magma allows to formalize the concept of Assembly Space in an alternative way to \cite{Marshall2022FormalisingTP}, where it is defined in terms of acyclic quivers (see Definition 11 in \cite{Marshall2022FormalisingTP}). The two approaches are related in a straightforward way via a mapping from the set of Assembly Addition Chains into the set of Acyclic Directed Graphs, such that the Assembly Index of an object $\mathscr{O}$ could be either defined as the length of any Optimal Assembly Addition Chain of $\mathscr{O}$ or as the length of any of the shortest Assembly Pathways of the object $\mathscr{O}$.         \\\\
This paper is organized as follows: In Section \ref{secAAC} we introduce the concept of Assembly Addition Chains and provide examples related to strings, graphs and polyominoes, including also a simple result regarding coarse lower and upper bounds for the length of Optimal Assembly Addition Chains. In Section  \ref{Secbound} we refine these bounds by deriving a better lower bound in terms of \textit{standard} OAC and better upper bounds based on Binary (see Definition \ref{DefBDO}) and 2-Piece (see Definition \ref{Def2PD}) decompositions of the objects. Then in Section \ref{secExa} we explicitly compute these upper bounds for the length of Optimal Assembly Addition Chains of Strings, Graphs and Polyominoes. Finally, we finish the paper presenting an outlook, conclusions and possible further research directions in Section \ref{secOut}. The explicit mapping between Assembly Addition Chains and Directed Acyclic Graphs is presented in Appendix \href{sec:AppA}{A}.\\\\
\textbf{Acknowledgments.} We are grateful to Anisha K. Yeddanapudi, Amit Kahana, Abhishek Sharma and Ian Seet for helpful discussions and motivating conversations. Special thanks to Stuart Marshall and Gage Siebert for carefully reading and providing useful comments about a preliminary draft of the present document. The work of L.C, J.C.M.P and K.Y.P is supported by the UKRI EPSRC, NASA, NIH, and the Templeton and the Sloan Foundations.

\section{Assembly Addition Chains}\label{secAAC}
Let $S$ be a set. As elucidated above, a generalized concept of Addition Chain over $S$ requires to define a \textit{gluing} binary operation over $S$. This cannot be simply a closed binary operation over $S$ since there could be different ways to \textit{glue} two objects. For example, if $S=2$-Strings (i.e. binary strings) then there are two ways to glue the directed strings 01010 and 001011, say $01010001011$ and $00101101010$. More drastically, if $S$ is the set of Connected Graphs there are multiple ways to glue two graphs by identifying different subsets of their vertices.\\\\
The algebraic structure that allows us to formally capture this concept of binary \textit{gluing} operation is an extension of the concept of Magma (see e.g. \cite{bergman}).
\begin{definition} \textbf{(Magma).} A \textit{Magma} $(S,\circ)$ consists of a set $S$ and a binary operation $\circ: S \times S \to S$. 
\end{definition}
Using this definition we introduce the algebraic concept of Multi-Magma. 
\begin{definition} \textbf{(Multi-Magma).} A \textit{Multi-Magma} $(S,\circ)$ consists of a set $S$ and a binary operation $\circ :2^S \times 2^S \to 2^S$ such that $(2^S,\circ)$ is a Magma.   
\end{definition}
Since our goal is to formally define Addition Chains over arbitrary sets $S$, we are particularly interested in Multi-Magma structures $(S,\circ)$ such that the binary operation $\circ$ takes two subsets $S_1$ and $S_2$ of $S$ and gives the subset $S_1 \circ S_2$ of $S$ that consists of \textit{all} the objects that could be \textit{constructed by gluing} together elements of $S_1$ with elements of $S_2$. More precisely, since for AC every step involves only the gluing of \textit{two} elements, like in (\ref{adi}), we are only interested in the restriction of the Multi-Magma structures over $S$ to the subset of $2^S$ that consists of the singletons of elements in $S$. \\\\
In the following we present Multi-Magma structures over the Sets of Integers, Strings, Graphs and Polyominoes. As remarked above, we explicitly show the binary operation for singletons, since is the only part of the structure required for defining AC over these sets; however, a  Multi-Magma structure over any of these sets is defined by extending it \textit{element-wise} for any pair of subsets $S_1$ and $S_2$.
\begin{exmp}
    $(S,\circ)$ such that $S=\mathbb{Z}^+$ and $\circ$ is defined by 
    \begin{equation}
        \{s\} \circ \{t\} \equiv \{s+t\}
    \end{equation}
for any $s,t \in \mathbb{Z}^+$. This Multi-Magma structure is just the element-wise sum of integers for any two subsets of positive integers.  
\end{exmp}

\begin{exmp}\label{exmpstrings}
Let $S=j$-Strings, i.e $S$ is the set of strings constructed using an alphabet of cardinality $j$. Over the set $S$ the group $\mathbb{Z}_2$ defines a natural action: if $s \in S$ then 
\begin{equation}
    0 \cdot s=s, \qquad  1 \cdot s=s^T
\end{equation}
where $s^T$ is the reversed string, e.g. for $j=2$ if $s=011101$ then $s^T=101110$.
 In terms of this action we define the sets of \textit{directed} $j$-Strings 
\begin{equation}
j\text{-Strings}_{\text{d}} \equiv j\text{-Strings} 
\end{equation}
and \textit{undirected} $j$-Strings 
\begin{equation}
 j\text{-Strings}_{\text{u}} \equiv j\text{-Strings}/\mathbb{Z}_2.
\end{equation}
In other words in the set of directed $j$-Strings every string is considered as a different object, meanwhile in the set of undirected $j$-Strings an string and its reverse are considered as the same object. \\\\
Over the set of $S=j$-Strings it is possible to define the semi-group structure given by the concatenation of strings, i.e if $s_1,s_2 \in j$-Strings then 
\begin{equation}
\begin{split}
    \frown: j\text{-Strings} \times j\text{-Strings} & \to j\text{-Strings} \\
    (s_1,s_2) & \mapsto s_1 \frown s_2=s_1s_2,
    \end{split}
\end{equation}
where \textit{juxtaposition} of $j$-Strings means that the two strings are glued together in the order presented, e.g. for $s_1=10101011$ and $s_2=101110$ then  
\begin{equation*}
   s_1 \frown s_2 = 10101011101110, \qquad s_2 \frown s_1= 10111010101011.
\end{equation*}
In terms of this semigroup structure we could define a Multi-Magma structure over $j\text{-Strings}_\text{d}$ given by
\begin{equation}
        \{s_1\} \circ \{s_2\} \equiv \{s_1 \frown s_2, s_2 \frown s_1\}
    \end{equation}
and a Multi-Magma structure over $j\text{-Strings}_{\text{u}}$ given by
\begin{equation}
\begin{split}
          \{s_1\} \circ \{s_2\} \equiv \{ \underbrace{s_1 \frown s_2}_{\sim s_2^T \frown s_1^T }, \underbrace{s_1^T \frown s_2}_{\sim s_2^T \frown s_1}, \underbrace{s_1 \frown s_2^T}_{ \sim s_2 \frown s_1^T}, \underbrace{s_1^T \frown s_2^T}_{\sim s_2 \frown s_1}\}.  
\end{split}
\end{equation}

As an example, in the case of directed $2$-Strings we have 
\begin{equation*}
    \{0010010\} \circ \{011101101\} = \{0010010011101101,0111011010010010 \}
\end{equation*}
and in the case of undirected $2$-Strings
\begin{equation*}
    \{01\} \circ \{110\}=\{01110, \underbrace{10110}_{\sim 01101}, \underbrace{01011}_{\sim 11010}, \underbrace{10011}_{\sim 11001} \}.
\end{equation*}
\end{exmp}

\begin{exmp}\label{exmgraphs}
Let $S=\text{CCG}$ be the set of Colored Connected Graphs. Over $S \times S \times S$ we define the following relation: $(G_1,G_2,G)$ is said to be a \textit{good triple} if:
\begin{itemize}
    \item $G_1$ is a subgraph of $G$ and its subgraph-complement is $G_2$.
    \item $G_2$ is a subgraph of $G$ and its subgraph-complement is $G_1$.
\end{itemize}
In terms of this relation it is possible to define a Multi-Magma structure over $S$ given by 
    \begin{equation}\label{MMGraphs}
        \{G_1\} \circ \{G_2\} \equiv \{G \in S \ | \ (G_1,G_2,G) \text{ is good}\},
    \end{equation}
for every $G_1,G_2 \in S$.  \\
For example the Multi-Magma structure (\ref{MMGraphs}) for the elements
\begin{tikzpicture}
  \draw[red,thick] (1,0.5) -- (1,0);
  \draw (0.5,0) -- (1,0);
  \draw (1,0) -- (1.5,0);
  \foreach \x/\y in {0.5/0,1/0,1.5/0,1/0.5} {
    \filldraw[black] (\x,\y) circle (1.5pt);
  }
\end{tikzpicture} and 
\begin{tikzpicture}
\tikzstyle{blueedge}=[line width=1pt, lightblue]
  \coordinate (A) at (0,0);
  \coordinate (B) at (1,0);
  \coordinate (C) at (0.5,0.75);
  \draw[blueedge] (A) -- (B);       
  \draw (B) -- (C);                   
  \draw (C) -- (A);                   
  \foreach \pt in {A,B,C} {
    \filldraw[black] (\pt) circle (1.5pt);
  }
\end{tikzpicture} in CSCPG is given by Figure \ref{fig:sumgraph}.

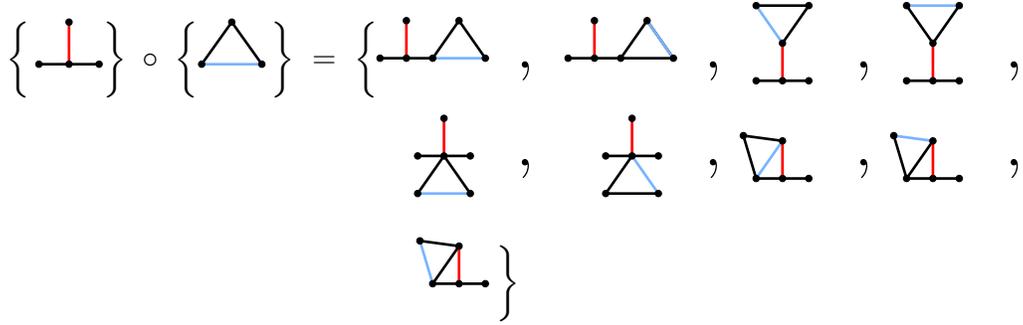
\begin{figure}[h]
    \centering
\begin{tikzpicture}[scale=1, every node/.style={inner sep=1pt}]
\tikzstyle{rededge}=[line width=1pt, red]
\tikzstyle{blueedge}=[line width=1pt, lightblue]
\tikzstyle{blackedge}=[line width=1pt]

\node at (-1,0) {$\bigg\{\, 
\begin{tikzpicture}[scale=0.8]
  \coordinate (a) at (0,0);
  \coordinate (b) at (1,0);
  \coordinate (d) at (0.5,0);
  \coordinate (c) at (0.5,0.7);
  \draw[blackedge] (a)--(d);
  \draw[blackedge] (d)--(b);
  \draw[rededge] (d)--(c);
    \foreach \pt in {a,b,c,d} {
    \filldraw[black] (\pt) circle (1.5pt);
  }
\end{tikzpicture}
\bigg\}\; \circ \; \bigg\{
\begin{tikzpicture}[scale=0.8]
  \coordinate (a) at (0,0);
  \coordinate (b) at (0.5,0.7);
  \coordinate (c) at (1,0);
  \draw[blueedge] (a)--(c);
  \draw[blackedge] (a)--(b)--(c);
    \foreach \pt in {a,b,c} {
    \filldraw[black] (\pt) circle (1.5pt);
  }
\end{tikzpicture}
\,\bigg\}  \;=\; \bigg\{ \rule{0pt}{12pt}$};


\begin{scope}[shift={(1.5,0)}]
  \coordinate (a) at (0,0);
  \coordinate (b) at (0.7,0);
  \coordinate (d) at (0.35,0);
  \coordinate (c) at (0.35,0.5);
  \coordinate (e) at (1.05,0.5);
  \coordinate (f) at (1.4,0);
  \draw[blueedge] (b)--(f);
  \draw[blackedge] (b)--(e)--(f);
  \draw[blackedge] (a)--(d);
  \draw[blackedge] (d)--(b);
  \draw[rededge] (d)--(c);
    \foreach \pt in {a,b,c,d,e,f} {
    \filldraw[black] (\pt) circle (1.2pt);
  }
\end{scope}
\node at (3.5,0) {\Huge$, \rule{0pt}{12pt}$};
\begin{scope}[shift={(4,0)}]
  \coordinate (a) at (0,0);
  \coordinate (b) at (0.7,0);
  \coordinate (d) at (0.35,0);
  \coordinate (c) at (0.35,0.5);
  \coordinate (e) at (1.05,0.5);
  \coordinate (f) at (1.4,0);
  \draw[blackedge] (b)--(f);
  \draw[blackedge] (b)--(e)--(f);
  \draw[blueedge] (e)--(f);
  \draw[blackedge] (a)--(d);
  \draw[blackedge] (d)--(b);
  \draw[rededge] (d)--(c);
    \foreach \pt in {a,b,c,d,e,f} {
    \filldraw[black] (\pt) circle (1.2pt);
  }
\end{scope}
\node at (6,0) {\Huge$, \rule{0pt}{12pt}$};
\begin{scope}[shift={(6.5,-0.3)}]
  \coordinate (a) at (0,0);
  \coordinate (b) at (0.7,0);
  \coordinate (d) at (0.35,0);
  \coordinate (c) at (0.35,0.5);
  \coordinate (e) at (0,1.0);
  \coordinate (f) at (0.7,1.0);
  \draw[blueedge] (c)--(e);
  \draw[blackedge] (e)--(f)--(c);
  \draw[blackedge] (a)--(d);
  \draw[blackedge] (d)--(b);
  \draw[rededge] (d)--(c);
    \foreach \pt in {a,b,c,d,e,f} {
    \filldraw[black] (\pt) circle (1.2pt);
  }
\end{scope}

\node at (8,0) {\Huge$, \rule{0pt}{12pt}$};
\begin{scope}[shift={(8.5,-0.3)}]
  \coordinate (a) at (0,0);
  \coordinate (b) at (0.7,0);
  \coordinate (d) at (0.35,0);
  \coordinate (c) at (0.35,0.5);
  \coordinate (e) at (0,1.0);
  \coordinate (f) at (0.7,1.0);
  \draw[blueedge] (e)--(f);
  \draw[blackedge] (f)--(c)--(e);
  \draw[blackedge] (a)--(d);
  \draw[blackedge] (d)--(b);
  \draw[rededge] (d)--(c);
    \foreach \pt in {a,b,c,d,e,f} {
    \filldraw[black] (\pt) circle (1.2pt);
  }
\end{scope}

\node at (10,0) {\Huge$, \rule{0pt}{12pt}$};


\begin{scope}[shift={(2,-1.3)}]
  \coordinate (a) at (0,0);
  \coordinate (b) at (0.7,0);
  \coordinate (d) at (0.35,0);
  \coordinate (c) at (0.35,0.5);
  \coordinate (e) at (0,-0.5);
  \coordinate (f) at (0.7,-0.5);
  \draw[blueedge] (e)--(f);
  \draw[blackedge] (f)--(d)--(e);
  \draw[blackedge] (a)--(d);
  \draw[blackedge] (d)--(b);
  \draw[rededge] (d)--(c);
    \foreach \pt in {a,b,c,d,e,f} {
    \filldraw[black] (\pt) circle (1.2pt);
  }
\end{scope}
\node at (3.5,-1.3) {\Huge$, \rule{0pt}{12pt}$};
\begin{scope}[shift={(4.5,-1.3)}]
  \coordinate (a) at (0,0);
  \coordinate (b) at (0.7,0);
  \coordinate (d) at (0.35,0);
  \coordinate (c) at (0.35,0.5);
  \coordinate (e) at (0,-0.5);
  \coordinate (f) at (0.7,-0.5);
  \draw[blueedge] (d)--(f);
  \draw[blackedge] (d)--(e)--(f);
  \draw[blackedge] (a)--(d);
  \draw[blackedge] (d)--(b);
  \draw[rededge] (d)--(c);
    \foreach \pt in {a,b,c,d,e,f} {
    \filldraw[black] (\pt) circle (1.2pt);
  }
\end{scope}
\node at (6,-1.3) {\Huge$, \rule{0pt}{12pt}$};
\begin{scope}[shift={(6.5,-1.6)}]
  \coordinate (a) at (0,0);
  \coordinate (b) at (0.7,0);
  \coordinate (d) at (0.35,0);
  \coordinate (c) at (0.35,0.5);
  \coordinate (e) at (-0.17,0.57);
  \draw[blueedge] (a)--(c);
  \draw[blackedge] (a)--(e)--(c);
  \draw[blackedge] (a)--(d);
  \draw[blackedge] (d)--(b);
  \draw[rededge] (d)--(c);
    \foreach \pt in {a,b,c,d,e} {
    \filldraw[black] (\pt) circle (1.2pt);
  }
\end{scope}
\node at (8,-1.3) {\Huge$, \rule{0pt}{12pt}$};
\begin{scope}[shift={(8.5,-1.6)}]
  \coordinate (a) at (0,0);
  \coordinate (b) at (0.7,0);
  \coordinate (d) at (0.35,0);
  \coordinate (c) at (0.35,0.5);
  \coordinate (e) at (-0.17,0.57);
  \draw[blueedge] (e)--(c);
  \draw[blackedge] (e)--(a)--(c);
  \draw[blackedge] (a)--(d);
  \draw[blackedge] (d)--(b);
  \draw[rededge] (d)--(c);
    \foreach \pt in {a,b,c,d,e} {
    \filldraw[black] (\pt) circle (1.2pt);
  }
\end{scope}
\node at (10,-1.3) {\Huge$, \rule{0pt}{12pt}$};



\begin{scope}[shift={(2.2,-3)}]
  \coordinate (a) at (0,0);
  \coordinate (b) at (0.7,0);
  \coordinate (d) at (0.35,0);
  \coordinate (c) at (0.35,0.5);
  \coordinate (e) at (-0.17,0.57);
  \draw[blueedge] (e)--(a);
  \draw[blackedge] (e)--(c)--(a);
  \draw[blackedge] (a)--(d);
  \draw[blackedge] (d)--(b);
  \draw[rededge] (d)--(c);
    \foreach \pt in {a,b,c,d,e} {
    \filldraw[black] (\pt) circle (1.2pt);
  }
\end{scope}


\node at (3.2,-3) {$\bigg\} \rule{0pt}{12pt}$};

\end{tikzpicture}
    \caption{Multi-Magma structure over the set of Colored Connected Graphs.}
    \label{fig:sumgraph}
\end{figure}
\end{exmp}

\begin{exmp}\label{CpolyMM}
The last example that we consider is Multi-Magma structures over the set $S=\text{CPoly}$ of (colored) Polyominoes, i.e. plane geometric figures formed by joining equal (colored) squares edge to edge.\footnote{More precisely in this paper by Polyominoes we mean \textit{free Polyominoes}, i.e. we identify two polyominoes if one could be obtained from the other by a rotation or/and a reflection, see e.g. \cite{golomb} and \cite{REDELMEIER1981191}.}  In Figure \ref{fig:polyexmp} we present all the Polyominoes of sizes three (trominoes), four (tetrominoes), five (pentominoes) and six (hexominoes). 

\begin{figure}[h]
    \centering
\includegraphics[width=0.8\textwidth]{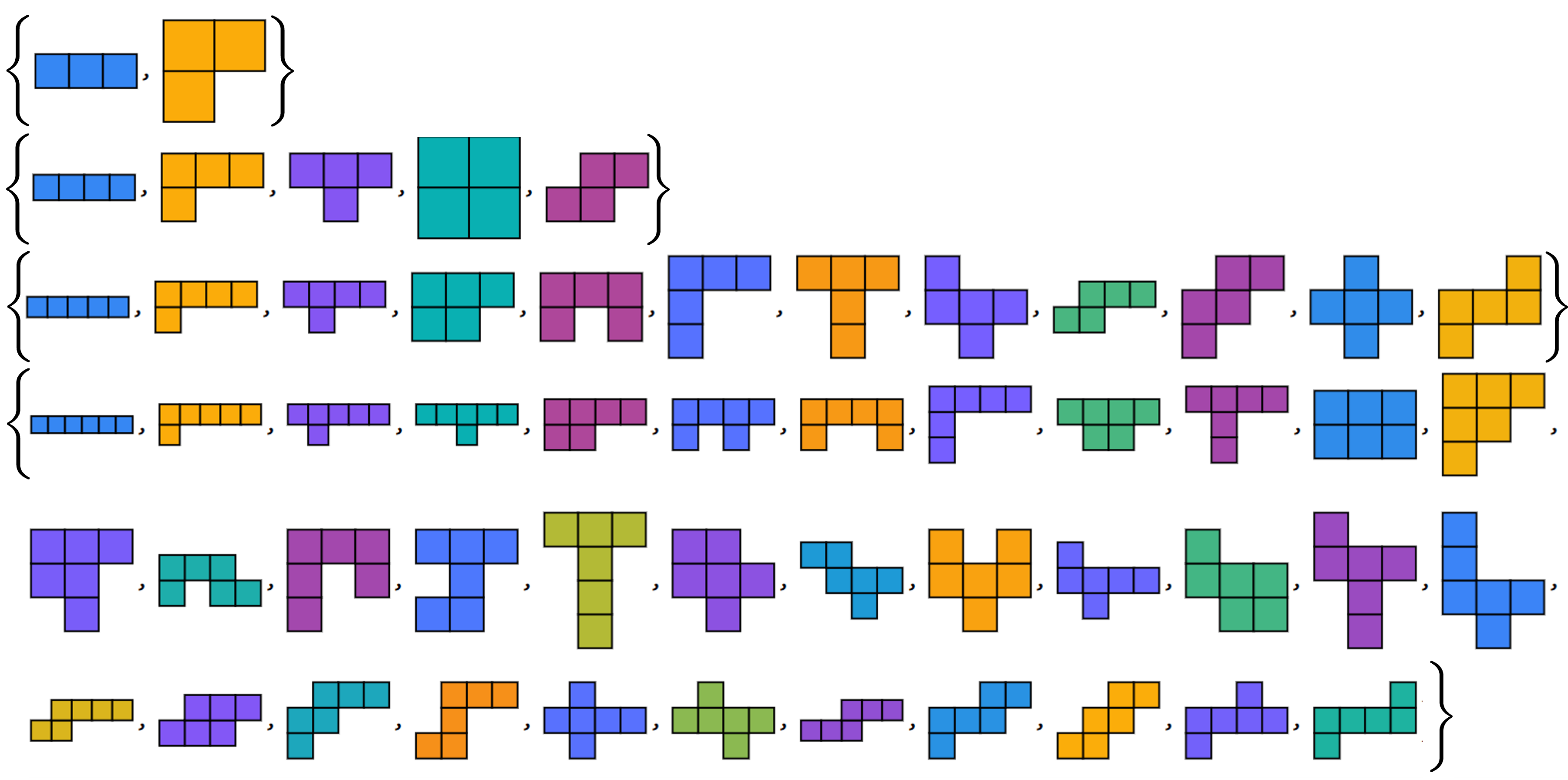}
    \caption{All Polyominoes of sizes three, four, five and six.}
    \label{fig:polyexmp}
\end{figure}

Analogously to the previous example, over $\text{Poly}\times\text{Poly}\times\text{Poly}$ we define the following relation: $(P_1,P_2,P)$ is said to be a \textit{good triple} if:
\begin{itemize}
    \item $P_1$ is a sub-Polyomino of $P$, and its Polyomino-complement is $P_2$.
    \item $P_2$ is a sub-Polyomino of $P$, and its Polyomino-complement is $P_1$.
\end{itemize}
In terms of this relation it is possible to define a Multi-Magma structure over $S$ given by 
    \begin{equation}\label{multipoly}
        \{P_1\} \circ \{P_2\} \equiv \{P \in S \ | \ (P_1,P_2,P) \text{ is good}\},
    \end{equation}
for every $P_1,P_2 \in S$. \\\\ 
For example the Multi-Magma structure (\ref{multipoly}) for the polyominoes \begin{tikzpicture}
  \foreach \x/\y in {0/0,0/0.3125,0/0.625,0.3125/0} {
    \draw[thick] (\x,\y) rectangle ++(0.3125,0.3125);
  }
\end{tikzpicture} and \begin{tikzpicture}
  \foreach \x/\y in {0/0,0/0.3125,0.3125/0,0.3125/0.3125} {
    \draw[thick] (\x,\y) rectangle ++(0.3125,0.3125);
  }
\end{tikzpicture} 
is given by Figure \ref{fig:sumpoly}.
\begin{figure}[h]
    \centering
\includegraphics[width=0.8\textwidth]{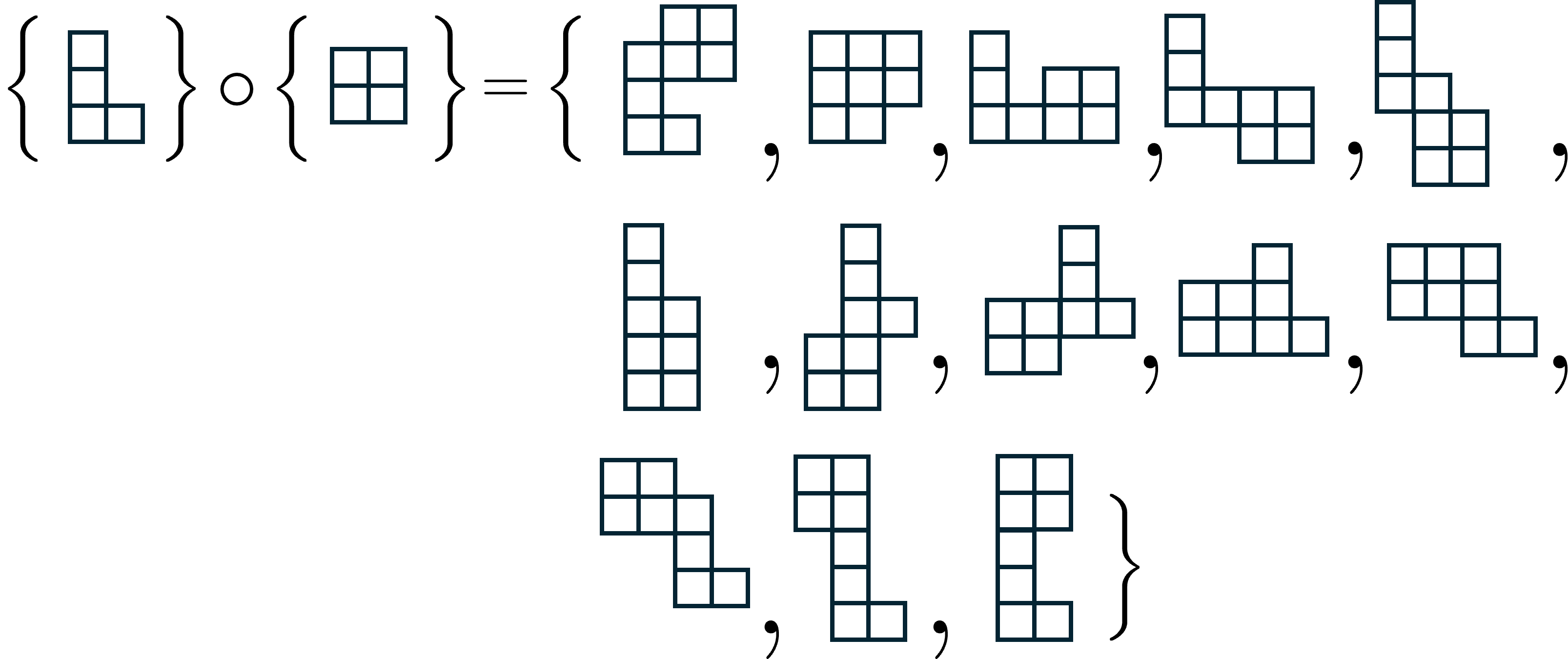}
    \caption{Multi-Magma structure over the set of Colored Fixed Polyominoes.}
    \label{fig:sumpoly}
\end{figure}
\end{exmp}

The algebraic structure of Multi-Magma $(S,\circ)$, analogously to the semi-group structure $(\mathbb{Z}^+,+)$, allows to define the concept of \textit{Assembly} Addition Chains over $S$. But meanwhile in the AC case there is a natural Building Block (say $1$) for a general Multi-Magma it is required to \textit{fix} a subset of $S$ as the set of Building Blocks for the \textit{Assembly} Addition Chains. 

\begin{definition}\label{AAC} \textbf{(Assembly Addition Chains).}
Let $(S,\circ)$ be a Multi-Magma, $BB \subset S$ (called the set of \textit{Building Blocks}) and $\mathscr{O} \in S \setminus BB$ an arbitrary object. Consider a sequence
\begin{equation*}
    (\mathscr{O}_1,\mathscr{O}_2,\dots, \mathscr{O}_r=\mathscr{O})
\end{equation*}
such that $\mathscr{O}_1 \in BB \circ BB$ and for every object $\mathscr{O}_\rho \in S$ for $\rho \in \{2,3,\dots,r\}$ there exist  $\mathscr{O}_\sigma,\mathscr{O}_\tau \in \{\mathscr{O}_1,\dots, \mathscr{O}_{\rho-1}\} \cup BB$ such that 
\begin{equation}
    \mathscr{O}_\rho \in \{\mathscr{O}_\sigma\} \circ \{\mathscr{O}_\tau\}.
\end{equation}
Such sequences will be called \textit{Assembly Addition Chains} (AAC) of $\mathscr{O} \in S$ with seeds on $BB$ and the set of all such sequences will be denoted by $AAC(\mathscr{O},BB)$. 
\end{definition}

This definition shows how the concept of Multi-Magma naturally allows to extend the notion of AC over general sets. Nevertheless, there is a particular property that is implicit in the semi-group $(\mathbb{Z}^+,+)$ but is not necessarily satisfied over an arbitrary Multi-Magma: any $n \in \mathbb{Z}^+$ could be represented by  an AC, say $(1,2,\dots,n)$, i.e. that $n$ can be represented as the sum of $n$ 1's. We define a family of Multi-Magma structures where an analogous property holds, i.e. a structure for which is guaranteed that \textbf{every} object is \textit{composed} by a finite number of Building Blocks. 

\begin{definition}\textbf{(Assembly Multi-Magma).}
Let $(S,\circ)$ be a Multi-Magma and $BB \subset S$. If for every $\mathscr{O} \in S \setminus BB$ there exists $\mathscr{B}_i \in BB$ for $i=0,1,\dots,r$ such that:
\begin{enumerate}
    \item 
    \begin{equation}\label{AMMC}
        \Gamma(\mathscr{O}) \equiv (\mathscr{O}_1,\mathscr{O}_2, \cdots,\mathscr{O}_r=\mathscr{O}) \in AAC(\mathscr{O},BB)
    \end{equation}
    with
    \begin{itemize}
        \item $  \mathscr{O}_1 \in \{ \mathscr{B}_0 \} \circ \{ \mathscr{B}_1 \}$,
        \item $\mathscr{O}_i \in \{\mathscr{O}_{i-1}\} \circ \{\mathscr{B}_i\}$ for $i=2,\dots,r$.
    \end{itemize}
\item If there exist $\overline{\mathscr{B}}_j \in BB$ for $j=0,1,\dots,s$ such that 
\begin{equation*}
\overline{\Gamma}(\mathscr{O})=(\overline{\mathscr{O}}_1,\overline{\mathscr{O}}_2,\cdots,\overline{\mathscr{O}}_s=\mathscr{O}) \in AAC(\mathscr{O},BB)
\end{equation*} 
with 
    \begin{itemize}
        \item $  \overline{\mathscr{O}}_1 \in \{ \overline{\mathscr{B}}_0 \} \circ \{ \overline{\mathscr{B}}_1 \}$,
        \item $\overline{\mathscr{O}}_i \in \{\overline{\mathscr{O}}_{i-1}\} \circ \{\overline{\mathscr{B}}_i\}$ for $i=2,\dots,r$,
    \end{itemize}
    then $r=s$ and there exists $\sigma \in S_{r}$ (the symmetric group of degree $r$) such that $\overline{\mathscr{B}}_i=\mathscr{B}_{\sigma(i)}$ for $i=0,1,\dots,r$, 
\end{enumerate}
then we say that the triple $(S,\circ,BB)$ is an \textit{Assembly Multi-Magma} (AMM).   
\end{definition}
Intuitively, if we interpret the binary operation $\circ$ as a \textit{gluing} operation, an AMM structure $(S,\circ,BB)$ has the property that every object $\mathscr{O} \in S$ is either (i) a Building Block or (ii) could be constructed by recursively \textit{gluing} together (via $\circ$) a \textit{unique set} of Building Blocks. In other words, an AMM structure over $S$ implies that every object $\mathscr{O}$ could be represented by an unique finite number of Building Blocks, such that the object could be \textit{constructed} by assembling recursively these Building Blocks. For any AMM $(S,\circ,BB)$ we denote by $\Gamma(\mathscr{O})$ any characteristic AAC of $\mathscr{O} \in S$ given by (\ref{AMMC}).
\begin{remark}
If the set of Building Blocks of an AMM $(S,\circ,BB)$ is known by context (as is the case in the examples considered here) we simply denote the set $AAC(\mathscr{O},BB)$ by $AAC(\mathscr{O})$ for every $\mathscr{O} \in S$. For example: 
\begin{enumerate}
    \item For $S=\mathbb{Z}^+$, then $BB=\{1\}$.
    \item For $S=j$-Strings, then $BB=\{0,1,2,...,j-1\}$.
    \item For $S=\text{CCG}$, then $BB$ is the set of monochromatic segments. 
    \item For $S=\text{CPoly}$, then $BB$ is the set of monochromatic $1 \times 1$ squares.
\end{enumerate}
\end{remark}
Analogously to the case of $AC(n)$ in $(\mathbb{Z}^+,+)$, over $AAC(\mathscr{O},BB)$ for any Multi-Magma $(s,\circ,BB)$ we could define the length of an Assembly Addition Chain. For every $\mathscr{O} \in S \setminus BB$, we define the \textit{length function}  
\begin{equation*}
\begin{split}
L: AAC(\mathscr{O},BB) & \to \mathbb{Z}^{+} \\
    (\mathscr{O}_1,\mathscr{O}_2,\dots, \mathscr{O}_r=\mathscr{O}) & \mapsto r
\end{split}
\end{equation*}
which could be naturally extended to all
\begin{equation*}
    AAC(S,BB) \equiv \cup_{\mathscr{O} \in S \setminus BB}AAC(\mathscr{\mathscr{O}},BB).
\end{equation*}
This function simply corresponds to the length of the sequences in $AAC(S)$. Heuristically, since the sequences in $AAC(\mathscr{O},BB)$ represent paths that indicate how to recursively construct $\mathscr{O}$ using elements of $BB$, the function $L$ indicates the number of steps of each of these paths. Particularly, if $(S,\circ,BB)$ is an AMM then if we evaluate the function $L$ on $\Gamma(\mathscr{O}) \in AAC(\mathscr{O},BB)$ we get a characteristic number of the object $\mathscr{O} \in S$, say its \textit{size}.
\begin{definition} \textbf{(Size of an object).} Let $(S,\circ,BB)$ be an AMM. The \textit{size function} over $S$ is defined by 
\begin{equation}\label{sizefun}
\begin{split}
    s: S & \to \mathbb{Z}^+ \\
    \mathscr{O} & \mapsto \begin{cases}
        1 & \quad \text{for} \quad \mathscr{O} \in BB, \\
        L(\Gamma(\mathscr{O}))+1 & \quad \text{for} \quad \mathscr{O} \in S \setminus BB. 
    \end{cases}
\end{split}
\end{equation} 
\end{definition}
As mentioned in the introduction, the main motivation behind the concepts of Multi-Magma and Assembly Multi-Magma is to formalize the notion of \textit{gluing} objects and \textit{assembling} objects by recursive gluing starting from Building Blocks, respectively. Given that the size function allows to determine the number of Building Blocks required to construct any object $\mathscr{O}$, there must be a compatibility condition between the Multi-Magma binary operation $\circ$ and the function $s$, since the size of an object constructed by gluing two objects should be the sum of the sizes of the later.            
\begin{definition} \textbf{(Assembly Space).} An Assembly Space $(S,\circ,BB)$ is an AMM such that the size function satisfies the \textit{additivity condition}, i.e. if $\mathscr{O} \in S \setminus BB$ and $\mathscr{O} \in \{\mathscr{O}_\sigma\} \circ \{\mathscr{O}_\tau \}$ then 
\begin{equation}\label{sizeprop}
    s(\mathscr{O})= s(\mathscr{O}_\sigma)+s(\mathscr{O}_\tau).
\end{equation}
\end{definition} 
Analogously to the function $\ell$ for the case of AC in $(\mathbb{Z}^+,+)$, we define the \textit{Assembly Index} function $a$ for the case of AAC in any Assembly Space $(S,\circ,BB)$. 
\begin{definition}\label{AI} \textbf{(Assembly Index).} Let $(S,\circ,BB)$ be an Assembly Space and $\mathscr{O} \in S$. The \textit{Assembly Index} (AI) of $\mathscr{O}$ is defined by  
    \begin{equation}\label{AI}
        a(\mathscr{O}) \equiv \begin{cases}
        0 \qquad & \text{for} \quad \mathscr{O} \in BB \\
        \min (L(AAC(\mathscr{O}))) \quad & \text{for} \quad \mathscr{O} \in S \setminus BB. \\
        \end{cases}      
    \end{equation}
\end{definition}

The Assembly Addition Chains in an Assembly Space are the mathematical objects that capture the intuitive notion of gluing in a recursive way Building Blocks in order to produce an object. The concept of Assembly Space was already introduced in \cite{Sharma2023AssemblyTE} with the main goal of describing the space of all the molecules that could be constructed in a recursive way from a specific set of bonds (the so-called Chemical Combinatorial Space). Indeed the properties that define an \textit{object} in \cite{Sharma2023AssemblyTE} matches with the properties that an object $\mathscr{O} \in S$ in an Assembly Space $(S,\circ,BB)$ exhibits: An object is \textit{finite} (given by $s(\mathscr{O})$), is \textit{distinguishable} and is \textit{breakable} such that the set of constraints to construct it from elementary building blocks is quantifiable (given by $a(\mathscr{O})$). \\\\
Finally, Optimal Assembly Addition Chains are defined exactly in the same way as Optimal Addition Chains are defined, just with the function $a$ replacing $\ell$ . 
\begin{definition} \textbf{(Optimal Assembly Addition Chains).} Let $(S,\circ,BB)$ be an Assembly Space and $\mathscr{O} \in S \setminus BB$. An Assembly Addition Chain $C \in AAC(\mathscr{O})$ is \textit{optimal} if 
\begin{equation}
    L(C)= a(\mathscr{O}).
\end{equation}
\end{definition}
    
\begin{remark}\label{impremark}
In \cite{Marshall2022FormalisingTP}, together with the concept of Assembly Spaces the notion of Assembly Pathway was originally introduced. As described mathematically in \cite{Marshall2022FormalisingTP}, given a (sub)set of Building Blocks $BB$ and an object $\mathscr{O}$ composed by all the elements of $BB$, an Assembly Pathway of $\mathscr{O}$ is a Directed Acyclic Graphs with $a(\mathscr{O})+\text{Card}(BB)$ vertices such that: (I) \textit{exactly} $\text{Card}(BB)$ vertices have indegree $0$ (the vertices corresponding to the Building Blocks) while \textit{all} the other $a(\mathscr{O})$ vertices have indegree two (the vertices corresponding to the objects constructed recursively in order to assemble $\mathscr{O}$) and (II) \textit{only one} vertex has outdegree $0$ (the vertex corresponding to $\mathscr{O}$). $AAC(\mathscr{O},BB)$ provide an algebraic presentation of the Assembly Paths of $\mathscr{O} \in S$ in an Assembly Space $(S,\circ,BB)$ and in Appendix \href{AppA}{A} an explicit mapping (projection) formalizing this relation is presented.
\end{remark}

The concepts of (i) Assembly Addition Chain $AAC(\mathscr{O})$, (ii) size $s(\mathscr{O})$ and (iii) Assembly Index $a(\mathscr{O})$ of any object $\mathscr{O} \in S$ in an Assembly Space, have been defined as extensions of the concepts of (i) Addition Chain $AC(n)$, (ii) size $n$ and (iii) $\ell(n)$ of any $n \in \mathbb{Z}^+$, respectively. Hence, we expect that for the former case there should be an inequality of the type (\ref{ACineq}), bounding the values of $a(\mathscr{O})$ in terms of $s(\mathscr{O})$. Following the discussion in Section 4 of the Mathematical Supplement of \cite{Sharma2023AssemblyTE}, it is straightforward to derive an inequality that provides coarse estimates of $a(\mathscr{O})$.

\begin{proposition}\label{prop1} Let $(S,\circ,BB)$ be an Assembly Space. The Assembly Index and the size of any objects $\mathscr{O} \in S$ satisfy the following inequalities
    \begin{equation}\label{ineq1}
        \log_2(s(\mathscr{O})) \leq a(\mathscr{O}) \leq s(\mathscr{O})-1.
    \end{equation}
\end{proposition}

\begin{proof} If $\mathscr{O} \in BB$, then the inequality holds trivially since $a(\mathscr{O})=0$ and $s(\mathscr{O})=1$. If $\mathscr{O} \in S \setminus BB$, then since $\Gamma( \mathscr{O}) \in AAC(\mathscr{O})$
we have 
\begin{equation*}
  a(\mathscr{O}) \equiv \min (L(AAC(\mathscr{O})))  \leq L(\Gamma(\mathscr{O})) = s(\mathscr{O})-1,  
\end{equation*}
getting the upper bound in (\ref{ineq1}). On the other hand, let $(\mathscr{O}_1,\mathscr{O}_2,\cdots,\mathscr{O}_a\equiv \mathscr{O}) \in    O   AAC(\mathscr{O})$. As a consequence of the additivity condition of the size function we have 
\begin{equation*}
    s(\mathscr{O}) \leq 2 \max_{i=1,...,a-1}\{ s(\mathscr{O}_i) \} \leq 2^{a(\mathscr{O})},
\end{equation*}
obtaining the lower bound in (\ref{ineq1}).
\end{proof}

In the case of standard Addition Chains, a generalization of the problem of computing the number $\ell(n)$ is well-known to be NP-complete, see e.g. \cite{NPAC}, \cite{Lehman}, \cite{LehmanPhD}. Currently there is not known exact solution to the problem, even though some algorithms (e.g. the so-called M-ary method) run in time $polylog(n)$ and give a $1+ o(1)$-approximation. Recently this result has been extended to the computation of the AI, proving that given an Assembly Space $(S,\circ,BB)$ and $\mathscr{O} \in S$ the problem of computing $a(\mathscr{O})$ is an NP-complete problem (see the Supplementary Information section of \cite{kempes2024assembly}). \\\\
Similarly to the case of Addition Chains, there is \textbf{not} an exact solution to the problem of computing the AI of an arbitrary object, even though several algorithms have been developed and implemented (see e.g.    \cite{seet2024rapid}). Nevertheless, up to now, any implementation of the algorithms to compute the Assembly Index in particular Assembly Spaces, e.g. Strings and CCG exhibit important limitations, particularly to compute the AI of \textit{big} objects (e.g. $S \in 2$-Strings such that $s(S)>32$ and $G \in \text{CCG}$ such that $s(G)>20$). Hence, having inequalities like (\ref{ACineq}) for $\ell$ and (\ref{ineq1}) for $a$ are of computational relevance, since at least these provide (rough) estimates. 

\section{Bounds for the Assembly Index}\label{Secbound}
\begin{figure}[h!]
    \centering
    \begin{minipage}{0.48\textwidth}
        \centering
        \includegraphics[width=\textwidth]{StringsFull.png}
        \caption{Assembly Index as a function of the string length for all the 2-Strings of length less or equal than 25.}
        \label{fig:allstrings}
    \end{minipage}\hfill
    \begin{minipage}{0.48\textwidth}
        \centering
        \includegraphics[width=\textwidth]{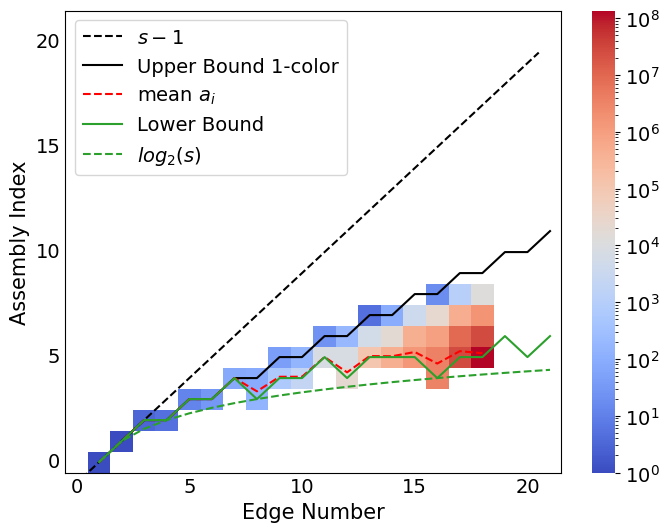}
        \caption{Assembly Index as a function of the edge number for all monochromatic SCG of edges less or equal than 18.}
        \label{fig:allgraphs}
    \end{minipage}
\end{figure}

In Figure \ref{fig:allstrings} and \ref{fig:allgraphs}, the AI of all the $\mathscr{O} \in  2$-Strings up to length 25 and all the $\mathscr{O} \in  \text{SCG}$ (Simple Connected  Graphs) up to length 18 are presented. The gaps between the bounds $s(\mathscr{O})-1$ and $\log_2(s(\mathscr{O}))$ and the maximum and minimum values of the AI, respectively, show that even though the bounds in (\ref{ineq1}) are correct they do not provide optimal estimates (specially for \textit{big} objects).

In this section we present two constructions that provide better upper and lower bounds for the values of the AI function $a$ (\ref{AI}) for any Assembly Space $(S,\circ,BB)$ in terms of the size function $s$ (\ref{sizefun}).

\subsection{Universal Lower Bound}
Due to the additivity property (\ref{sizeprop}) of the size function over an Assembly Space $(S,\circ,BB)$ and $s|_{BB}\equiv 1$, there exists a map 
\begin{equation}\label{mapF}
\begin{split}
    F: AAC(\mathscr{O},BB) & \to AC(s(\mathscr{O})) \\(\mathscr{O}_1,\mathscr{O}_2,\dots,\mathscr{O}_r=\mathscr{O}) & \mapsto \text{ord}[(s(\mathscr{O}_1),s(\mathscr{O}_2),\dots,s(\mathscr{O}))]
\end{split}
\end{equation}
where $\text{ord}$ is the function that orders the positive integers in increasing order and removes possible repetitions. Consequently, we have that for any $\mathscr{O} \in S$ and $C \in AAC(\mathscr{O},BB)$, the map $F$ satisfies the monotonicity condition 
\begin{equation}
    L(F(C)) \leq L(C).
\end{equation}
This injection of $AAC(\mathscr{O},BB)$ into $AC(s(\mathscr{O}))$, or more generally of $\cup_{\mathscr{O}\in S}AAC(\mathscr{O},BB)$ into $\cup_{n \in \mathbb{Z}^+}AC(n)$, implies that the AI of any object $\mathscr{O} \in S$ is bounded below by $\ell(s(\mathscr{O}))$. 

\begin{theorem}\label{Theo1} Let $(S,\circ,BB)$ be an Assembly Space. For any $\mathscr{O} \in S$ the following inequalities hold
    \begin{equation}\label{ineqbadai}
        \ceil*{\log_2 (s(\mathscr{O}))} \leq \ell(s(\mathscr{O})) \leq a(\mathscr{O}).
    \end{equation}
\end{theorem}

\begin{figure}[h]
    \centering
\includegraphics[width=0.6\textwidth]{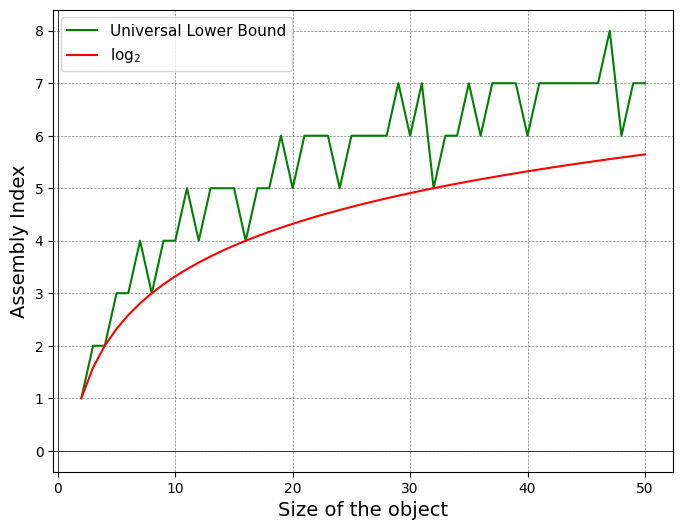}
    \caption{Plot of $\ell(s(\mathscr{O}))$ (green curve) and $\log_2(s(\mathscr{O}))$ (red curve) for objects of sizes up to 50.}
    \label{fig:Min}
\end{figure}

\begin{proof}
In \cite{Scholz} Scholz proved that if $2^{m}+1 \leq n \leq 2^{m+1}$, then 
\begin{equation}
 m + 1 \leq \ell(n) \leq 2m.    
\end{equation}
Hence if $\mathscr{O} \in S$ is an object such that $s(\mathscr{O})=2^k$ for some $k \in \mathbb{Z}^+$ then 
\begin{equation*}
   \ceil*{\log_2 (s(\mathscr{O}))} = \log_2 (s(\mathscr{O})) = k \leq \ell(s(\mathscr{O})).
\end{equation*}
Similarly, if $\mathscr{O} \in S$ is such that $2^{m}+1 \leq s(\mathscr{O}) < 2^{m+1}$ for some $k \in \mathbb{Z}^+$ then 
\begin{equation*}
   \ceil*{\log_2 (s(\mathscr{O}))} = m+1 \leq \ell(s(\mathscr{O})),
\end{equation*}
proving the first inequality. \\\\
On the other hand, by using the map $F$ (\ref{mapF}) we get that for any $S \in OAAC(\mathscr{O})$
\begin{equation*}
   \ell(s(\mathscr{O})) \leq L(F(S)) \leq  L(S)=a(\mathscr{O}),
\end{equation*}
proving in this way the last inequality.
\end{proof}
Theorem \ref{Theo1} shows that for \textbf{any} object $\mathscr{O}$ in an Assembly Space $(S,\circ,BB)$ the number $\ell(s(\mathscr{O}))$ is a better lower bound of $a(\mathscr{O})$ than the lower bound $\log_2(s(\mathscr{O}))$ of Proposition \ref{prop1}. In Figure \ref{fig:Min} we plot $\log_2(s(\mathscr{O}))$ and $\ell(s(\mathscr{O}))$ finding, as follows from (\ref{ineqbadai}), that the later is a better lower bound for $a(\mathscr{O})$ except in the cases when $s(\mathscr{O})=2^k$ for $k \in \mathbb{Z}^+$ (values for which both bounds coincide).

\begin{remark}
For the examples of Assembly Spaces in Section \ref{secAAC}, say $j$-Strings, CCG and CPoly, it is possible to explicitly describe \textit{a} family of objects such that $a(\mathscr{O})=\ell(s(\mathscr{O}))$: This family corresponds to the \textit{monochromatic string-like} objects. The adjective \textit{monochromatic} means that only one type of Building Block is used to construct the object and \textit{string-like} means that the object is constructed by gluing the Building Blocks in a ``\textit{linear way}''. 
\end{remark}


\begin{exmp}
    In the case where $S=j$-Strings, the family of monochromatic string-like objects consists of all the strings of the form 
    \begin{equation*}
        S_{b,n} \equiv \underbrace{bb\cdots b}_{n}
    \end{equation*}
    for $b \in \{0,1,\dots,j-1\}$ and $n \in \mathbb{Z}^+$. Indeed, 
    \begin{equation*}
a(S_{b,n})=\ell(s(S_{b,n}))=\ell(n).
    \end{equation*}
\end{exmp}

\begin{exmp} In the case where $S=$CCG, the family of monochromatic string-like objects consists of all the graphs of the form  
\begin{center}
\begin{tikzpicture}
  \node at (-1,0) {$LG_n \equiv$};
  \draw (0,0) -- (1,0);
  \filldraw[black] (0,0) circle (2pt);
  \filldraw[black] (1,0) circle (2pt);
  \node[below] at (0.5,0) {\small 1};
  \draw (1,0) -- (2,0);
  \filldraw[black] (2,0) circle (2pt);
  \node[below] at (1.5,0) {\small 2};
  \node at (2.75,0) {$\cdots$};
  \draw (3.5,0) -- (4.5,0);
  \filldraw[black] (3.5,0) circle (2pt);
  \filldraw[black] (4.5,0) circle (2pt);
  \node[below] at (4,0) {\small $n$};
\end{tikzpicture}.
\end{center}
Clearly in this case
\begin{equation*}
    a(LG_{n})=\ell(s(LG_{n}))=\ell(n).
\end{equation*}
\end{exmp}

\begin{exmp} In the case where $S=$CPoly, the family of monochromatic string-like objects consists of all the Polyominoes of the form 
\begin{center}
\begin{tikzpicture}
  \node at (-1,0.375) {$LP_n \equiv$};
  \draw (0,0) rectangle (0.75,0.75);
  \node[below] at (0.375,0) {\small 1};
  \draw (0.75,0) rectangle (1.5,0.75);
  \node[below] at (1.125,0) {\small 2};
  \node at (2.1,0.375) {$\cdots$};
  \draw (2.5,0) rectangle (3.25,0.75);
  \node[below] at (2.875,0) {\small $n$};
\end{tikzpicture}
\end{center}
and for this type of objects we have
\begin{equation*}
    a(LP_{n})=\ell(s(LP_{n}))=\ell(n).
\end{equation*}
\end{exmp}


\subsection{Upper Bounds}

As mentioned above the problem of calculating the AI  of an arbitrary object $\mathscr{O} \in S$ in an Assembly space $(S,\circ,BB)$ is computationally expensive (see e.g. \cite{kempes2024assembly} and \cite{seet2024rapid}) and there is yet not a generic algorithm that allows to compute it exactly for an arbitrary Set. The possibility of computing a better upper bound for $a(\mathscr{O})$ than $s(\mathscr{O})-1$ has been pointed qualitatively by other authors (see e.g \cite{Marshall2022FormalisingTP} and \cite{math12101600}), however their results remained conjectures (see also Figure 4 in \cite{RSPaper}). In this section, we derive better upper bounds for two particular types of objects of any Assembly Space (say Binary and 2-Piece Decomposable Objects, see below for the Definitions), getting in this way a better upper bound than in (\ref{ineq1}) for all the examples in Section \ref{secAAC}.
\begin{remark}
    We fix the following notation: 
    \begin{itemize}
        \item Let $(S,\circ,BB)$ be an Assembly Space. We introduce the notation
\begin{equation}
    S(k) \equiv \{\mathscr{O} \in S \ | \ s(\mathscr{O})=k\},
\end{equation}
i.e. by $S(k)$ we denote the set of elements of $S$ with size $k \in \mathbb{Z}^+ $.  
\item For any $k \in \mathbb{Z}^+$ we denote by $B(k)=\{n_1,n_2,\cdots,n_{H(k)}\} \subset \mathbb{Z}_{\geq 0}$ the set of non-negative integer numbers such that $n_1>n_2>\cdots>n_{H(k)}$ and $k=2^{n_1}+2^{n_2}+\cdots+2^{n_{H(k)}}$, i.e. $B(k)$ is set of powers of 2 in the binary decomposition of $k$ such that $H(k)$ is the Hamming weight of $k$.
    \end{itemize}
\end{remark}
First we start by considering a particular family of objects in any Assembly Space. 
\begin{definition}\label{DefBDO}
    \textbf{(Binary Decomposable Objects).} An object $\mathscr{O} \in S$ in an Assembly space $(S,\circ,BB)$ is said to be \textit{Binary Decomposable} if there exists 
    \begin{equation}
        \mathscr{B}^{(i)}_{j} \in BB \qquad \text{with } n_i \in B(s(\mathscr{O})) \quad  \text{and} \quad j \in {1,\dots,2^{n_i}}
    \end{equation} 
    such that 
    \begin{equation}
        \mathscr{O} \in (\mathscr{P}_{H(s(\mathscr{O}))} \circ \cdots \circ  (\mathscr{P}_3 \circ (\mathscr{P}_{2} \circ \mathscr{P}_{1} \underbrace{)) \cdots)}_{H(s(\mathscr{O}))-1}
    \end{equation}
    where
    \begin{equation}\label{pieces}
        \mathscr{P}_i \in \underbrace{\big( \cdots \big( \big(}_{i} \mathscr{B}^{(i)}_1 \circ  \mathscr{B}^{(i)}_2 \big) \circ \big(\mathscr{B}^{(i)}_3 \circ \mathscr{B}^{(i)}_4 \big) \big) \circ \cdots \circ \big( \big( \mathscr{B}^{(i)}_{2^{n_i}-3} \circ  \mathscr{B}^{(i)}_{2^{n_i}-2} \big) \circ \big(  \mathscr{B}^{(i)}_{2^{n_i}-1} \circ  \mathscr{B}^{(i)}_{2^{n_i}}\underbrace{\big) \big) \cdots \big)}_{i}.
    \end{equation}
\end{definition} 
We denote by $BD(S)$ the set of Binary Decomposable objects of an Assembly Space $(S,\circ,BB)$. Intuitively speaking a Binary Decomposable object is an object that could be recursively constructed by gluing together sub-objects whose sizes are powers of two.

\begin{exmp}
For $S=j$-Strings, $S=BD(S)$ since any string of size $a+b$ could be constructed via the juxtaposition of sub-strings of lengths $a$ (the first $a$ characters of the string) and $b$ (the last $b$ characters of the string).
\end{exmp}

\begin{exmp}\label{nonexCSCPG}
For $S=$CCG, we have $BD(S) \neq \emptyset$ but the inclusion $BD(S) \subset S$ is strict. In Figure \ref{figbdgraphs} (a) we show an example of an element in BD(S), while in (b) an object $\mathscr{O} \in S$ such that $\mathscr{O} \notin BD(S)$ is presented.
\begin{figure}[h]
\begin{subfigure}[h]{0.6\linewidth}
\begin{tikzpicture}[
    scale=1,
    every node/.style={inner sep=1pt},
    dot/.style={circle,fill=black,inner sep=1.2pt},
    rededge/.style={line width=1pt, red},
    blueedge/.style={line width=1pt, lightblue},
    blackedge/.style={line width=1pt}
]


\newcommand{\GraphA}{
  \begin{tikzpicture}[scale=0.35]
\node[dot] (A) at (0,0) {};        
\node[dot] (B) at (1,0) {};
\node[dot] (C) at (1,1) {};    
\node[dot] (D) at (2.0,2.0) {};    
\node[dot] (N) at (2.0,3.0) {}; 
\node[dot] (E) at (3,1) {};    
\node[dot] (F) at (3,0) {};      
\node[dot] (G) at (2.0,-1.5) {};   
\node[dot] (H) at (2,-0.5) {};   
\node[dot] (I) at (4,0) {};   
\node[dot] (J) at (4,-1.0) {};
\node[dot] (K) at (3,-1.5) {};
\node[dot] (L) at (3,-2.5) {};
\node[dot] (M) at (4,-2.5) {};


\draw[rededge] (A) -- (B);

\draw[blackedge] (B) -- (C);

\draw[rededge] (C) -- (D);

\draw[blackedge] (D) -- (E);

\draw[blueedge] (C) -- (E);

\draw[blueedge] (E) -- (F);

\draw[blackedge] (B) -- (H);

\draw[blackedge] (H) -- (G);

\draw[blackedge] (D) -- (N);

\draw[blackedge] (F) -- (I);

\draw[blueedge] (I) -- (J);

\draw[rededge] (J) -- (K);

\draw[blackedge] (K) -- (L);

\draw[blackedge] (L) -- (M);

\draw[rededge] (H) -- (F);
  \end{tikzpicture}
}

\newcommand{\Pfour}{
  \begin{tikzpicture}[scale=0.35]

\node[dot] (C) at (1,0) {};    
\node[dot] (E) at (3,0) {};    
\node[dot] (F) at (3,-1) {};      
\node[dot] (H) at (2,-1.5) {};   
\node[dot] (I) at (4,-1) {};   
\node[dot] (J) at (4,-2.0) {};
\node[dot] (K) at (3,-2.5) {};
\node[dot] (L) at (3,-3.5) {};
\node[dot] (M) at (4,-3.5) {};


\draw[blueedge] (C) -- (E);

\draw[blueedge] (E) -- (F);

\draw[blackedge] (F) -- (I);

\draw[blueedge] (I) -- (J);

\draw[rededge] (J) -- (K);

\draw[blackedge] (K) -- (L);

\draw[blackedge] (L) -- (M);

\draw[rededge] (H) -- (F);
  \end{tikzpicture}
}

\newcommand{\Pthree}{
  \begin{tikzpicture}[scale=0.35]
\node[dot] (B) at (1,0) {};  
\node[dot] (C) at (1,1) {};    
\node[dot] (E) at (3,1) {}; 
\node[dot] (D) at (2.0,2.0) {};    
\node[dot] (N) at (2.0,3.0) {}; 


\draw[blackedge] (B) -- (C);

\draw[rededge] (C) -- (D);

\draw[blackedge] (D) -- (E);

\draw[blackedge] (D) -- (N);
  \end{tikzpicture}
}

\newcommand{\Ptwo}{
  \begin{tikzpicture}[scale=0.35]
\node[dot] (A) at (0,0) {};        
\node[dot] (B) at (1,0) {};  

\node[dot] (H) at (2,-0.5) {};   


\draw[rededge] (A) -- (B);

\draw[blackedge] (B) -- (H);

  \end{tikzpicture}
}

\newcommand{\Pone}{
  \begin{tikzpicture}[scale=0.35]

\node[dot] (G) at (2.0,-1.5) {};   
\node[dot] (H) at (2,-0.5) {};   


\draw[blackedge] (H) -- (G);

  \end{tikzpicture}
}

\newcommand{\PthreeOne}{
  \begin{tikzpicture}[scale=0.35]
\node[dot] (B) at (1,0) {};  
\node[dot] (C) at (1,1) {};    

\draw[blackedge] (B) -- (C);

  \end{tikzpicture}
}

\newcommand{\PthreeTwo}{
  \begin{tikzpicture}[scale=0.35]
\node[dot] (C) at (1,1) {};    
\node[dot] (D) at (2.0,2.0) {};    

\draw[rededge] (C) -- (D);

  \end{tikzpicture}
}

\newcommand{\PthreeThree}{
  \begin{tikzpicture}[scale=0.35]
\node[dot] (E) at (3,1) {}; 
\node[dot] (D) at (2.0,2.0) {};    

\draw[blackedge] (D) -- (E);

  \end{tikzpicture}
}

\newcommand{\PthreeFour}{
  \begin{tikzpicture}[scale=0.35]
\node[dot] (D) at (2.0,2.0) {};    
\node[dot] (N) at (2.0,3.0) {}; 
\draw[blackedge] (D) -- (N);
  \end{tikzpicture}
}

\newcommand{\PfourOne}{
  \begin{tikzpicture}[scale=0.35]

\node[dot] (C) at (1,0) {};    
\node[dot] (E) at (3,0) {};    

\draw[blueedge] (C) -- (E);

  \end{tikzpicture}
}

\newcommand{\PfourTwo}{
  \begin{tikzpicture}[scale=0.35]

\node[dot] (E) at (3,0) {};    
\node[dot] (F) at (3,-1) {};      

\draw[blueedge] (E) -- (F);

  \end{tikzpicture}
}

\newcommand{\PfourThree}{
  \begin{tikzpicture}[scale=0.35]

\node[dot] (F) at (3,-1) {};      
\node[dot] (H) at (2,-1.5) {};   

\draw[rededge] (H) -- (F);
  \end{tikzpicture}
}

\newcommand{\PfourFour}{
  \begin{tikzpicture}[scale=0.35]

\node[dot] (F) at (3,-1) {};      
\node[dot] (I) at (4,-1) {};   

\draw[blackedge] (F) -- (I);

  \end{tikzpicture}
}

\newcommand{\PfourFive}{
  \begin{tikzpicture}[scale=0.35]

\node[dot] (I) at (4,-1) {};   
\node[dot] (J) at (4,-2.0) {};
\draw[blueedge] (I) -- (J);

  \end{tikzpicture}
}

\newcommand{\PfourSix}{
  \begin{tikzpicture}[scale=0.35]

\node[dot] (J) at (4,-2.0) {};
\node[dot] (K) at (3,-2.5) {};

\draw[rededge] (J) -- (K);

  \end{tikzpicture}
}
\newcommand{\PfourSeven}{
  \begin{tikzpicture}[scale=0.35]

\node[dot] (L) at (3,-3.5) {};
\node[dot] (M) at (4,-3.5) {};

\draw[blackedge] (L) -- (M);

  \end{tikzpicture}
}
\newcommand{\PfourEight}{
  \begin{tikzpicture}[scale=0.35]

\node[dot] (K) at (3,-2.5) {};
\node[dot] (L) at (3,-3.5) {};

\draw[blackedge] (K) -- (L);

  \end{tikzpicture}
}


\node at (0,0) {\GraphA};
\node[font=\large] at (5.1,0) {$ \in \left( \underbrace{\Pfour}_{\mathscr{P}_4} \circ \left( \underbrace{\Pthree}_{\mathscr{P}_3} \circ \left( \underbrace{\Ptwo}_{\mathscr{P}_2} \circ \underbrace{\Pone}_{\mathscr{P}_1} \right)\right)\right)$};


\node[font=\large] at (2.8,-2.2) {$\Pthree \in \left(\, \PthreeOne \circ \PthreeTwo \,\right)\circ \left(\, \PthreeFour\circ \PthreeThree \,\right)$};


\node[font=\large] at (5.2,-4.2) {$\Pfour\in \left( \left(\PfourOne \circ \PfourTwo\right) \circ\left(\PfourThree \circ \PfourFour\right)\circ\left(\PfourSix \circ \PfourFive \right)\circ \left(\PfourEight \circ   \PfourSeven\right)\right)$};

\end{tikzpicture}
\caption{Binary Decomposable Colored Connected Graph of size 15.}
\end{subfigure}
\hfill
\begin{subfigure}[h]{0.2\linewidth}
\begin{tikzpicture}[
    dot/.style={circle,fill=black,inner sep=1.2pt},
    blackedge/.style={line width=1pt},
    rededge/.style={line width=1pt,red},
    blueedge/.style={line width=1pt,lightblue}
]

\node[dot] (A) at (0.2,0) {};
\node[dot] (Z) at (0.8,0) {};
\node[dot] (B) at (1.4,0) {};
\node[dot] (C) at (2.0,0) {};
\node[dot] (D) at (2.0,0.6) {};
\node[dot] (E) at (2.0,1.2) {};
\node[dot] (W) at (2.0,1.8) {};
\node[dot] (F) at (2.6,0) {};
\node[dot] (G) at (2.6,-0.6) {};

\draw[blackedge] (A) -- (Z);

\draw[rededge] (Z) -- (B);

\draw[blueedge] (B) -- (C);

\draw[blueedge] (C) -- (D);

\draw[blackedge] (D) -- (E);

\draw[rededge] (E) -- (W);

\draw[blackedge] (C) -- (F);

\draw[blackedge] (F) -- (G);

\end{tikzpicture}
\caption{Colored Connected Graph of size 8 that is not Binary Decomposable.}
\end{subfigure}%
\caption{An (a) example and a (b) non-example of Binary Decomposable Colored Connected Graphs.}
\label{figbdgraphs}
\end{figure}
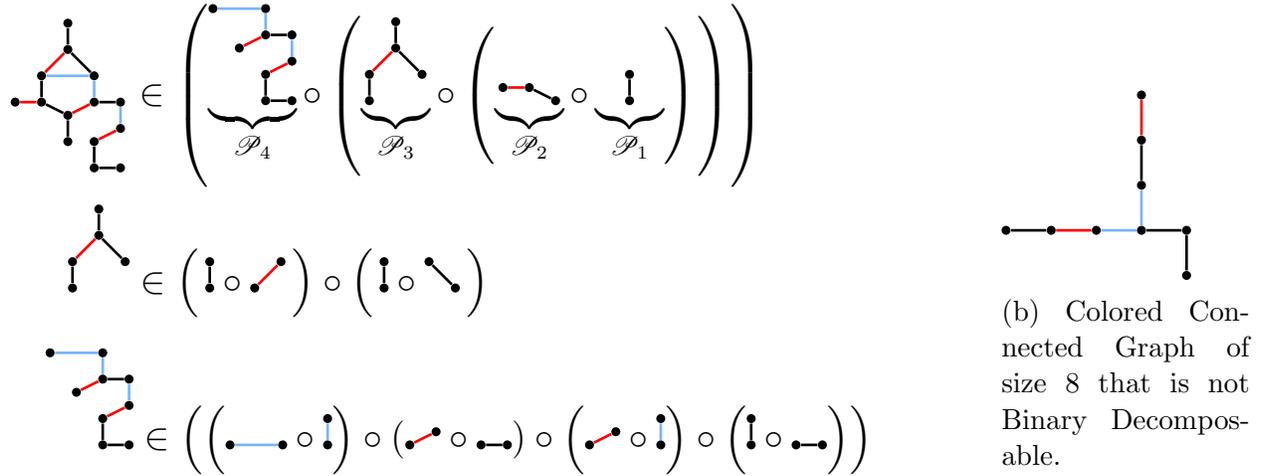
\end{exmp}

\begin{exmp}\label{nonexPoly}
Similarly to the previous case, for $S=$CPoly we have $BD(S)\neq \emptyset$ and $BD(S) \neq \emptyset$. In Figure \ref{figbdpoly} (a) we present $\mathscr{O} \in BD(S)$, while in (b) an object such that $\mathscr{O} \notin BD(S)$.
\begin{figure}[h]
\begin{subfigure}[h]{0.4\linewidth}
\centering
\includegraphics[width=0.5\linewidth]{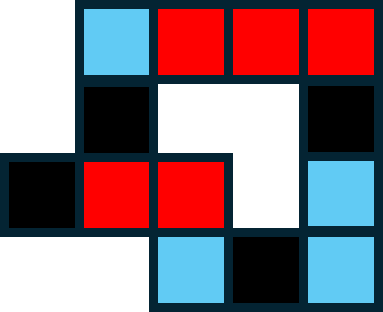}

\caption{Binary Decomposable Colored Polyomino of size 13.}
\end{subfigure}
\hfill
\begin{subfigure}[h]{0.4\linewidth}
\centering
\includegraphics[width=0.5\linewidth]{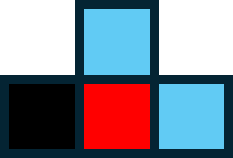}
\caption{Colored Polyomino of size 4 that is not Binary Decomposable.}
\end{subfigure}%
\caption{An (a) example and a (b) non-example of Binary Decomposable Colored Polyominoes.}
\label{figbdpoly}
\end{figure}
\end{exmp}

\begin{lemma}\label{lemmaBD}
    Let $(S,\circ,BB)$ be an Assembly Space and $\mathscr{O} \in BD(S)$ such that $s(\mathscr{O})=2^n$. Then 
    \begin{equation}\label{ineqBDObj}
        a(\mathscr{O}) \leq \sum^n_{i=1} \min \{2^{n-i},\text{Card}(S(2^i)) \}. 
    \end{equation}
\end{lemma}

\begin{proof}
  Since $\mathscr{O} \in BD(S)$ and $s(\mathscr{O})=2^n$, there exist $\mathscr{B}_i \in BB$ for $i=1,\dots,2^n$ such that 
\begin{equation*}
     \mathscr{O} \in \underbrace{\big( \cdots \big( \big(}_{n}\mathscr{B}_1 \circ \mathscr{B}_2 \big) \circ \big(\mathscr{B}_3 \circ \mathscr{B}_4 \big) \big) \circ \cdots \circ \big( \big(\mathscr{B}_{2^{n}-3} \circ \mathscr{B}_{2^{n}-2} \big) \circ \big( \mathscr{B}_{2^{n}-1} \circ \mathscr{B}_{2^{n}}\underbrace{\big) \big) \cdots \big)}_{n}.
\end{equation*} 
Hence there exists $C \in AAC(\mathscr{O})$ such that $L(C)=2^n-1$ explicitly given by
\begin{equation}\label{AACproof}
    C= (\mathscr{B}_1 \circ \mathscr{B}_2, \mathscr{B}_3 \circ \mathscr{B}_4,\cdots,\mathscr{B}_{2^n-1}\circ \mathscr{B}_{2^n},(\mathscr{B}_1 \circ \mathscr{B}_2)\circ(\mathscr{B}_3 \circ \mathscr{B}_4), \cdots, \mathscr{O}),
\end{equation}
i.e. where the first $2^{n-1}$ objects belong to $\mathscr{B}_i \circ \mathscr{B}_{i+1}$ for $i=1,3,\cdots,2^n-1$, respectively, and the next $2^{n-2}$ objects belong to $(\mathscr{B}_i \circ \mathscr{B}_{i+1}) \circ (\mathscr{B}_{i+2} \circ \mathscr{B}_{i+3})$ for $i=1,5,\cdots,2^{n}-3$, respectively, and so on until the last object is $\mathscr{O}$. Hence, the first $2^{n-1}$ objects in $C$ have size $2$ , the next $2^{n-2}$ objects have size $4$ and so on until the last object ($\mathscr{O}$) has size $2^n$. Since the maximum number of elements of size $k$ in $S$ is given by $\text{Card}(S(k))$, there exists $C' \in AAC(\mathscr{O})$ (obtained from $C$ by removing possible repeated objects) such that  
\begin{equation*}
    a(\mathscr{O}) \leq L(C') = \sum^n_{i=1} \min \{2^{n-i},\text{Card}(S(2^i)) \},
\end{equation*}
getting in this way an upper bound for the AI of any object $\mathscr{O} \in BD(S)$ whose size is a power of two.
\end{proof}
For the Assembly Spaces $S$ considered along this paper, i.e. $j$-Strings, CCG and CPoly, the cardinality of the level sets $S(k)$ increases monotonically with $k$. Hence for each $n \in \mathbb{Z}^+$ there exists a critical integer number  
\begin{equation}\label{c_n}
    c(n)= \min \{ \mu_n^{-1}(\mathbb{Z}_{\geq 0}) \},
\end{equation}
in terms of the function 
$\mu_n:\{0,1,2,\cdots,n-1\} \to \mathbb{Z}$ given by
\begin{equation}
    \mu_n(j)= 2^j - \text{Card}(S({2^{n-j}})),
\end{equation}
such that if $j > c(n)$ then  $2^{j}-\text{Card}(S(2^{n-j}))$ objects of size $2^{n-j}$ are repetitions in (\ref{AACproof}), and consequently for all these cases the inequality (\ref{ineqBDObj}) becomes
\begin{equation}\label{v2ineqBD}
    a(\mathscr{O}) \leq 2^{c(n)}-1 + \sum_{i=1}^{n-c(n)} \text{Card}(S(2^{i})).
\end{equation}
\begin{remark}\label{remrk6}
The counting corresponding to the RHS of inequality (\ref{v2ineqBD}) is depicted in Figure \ref{fig:MaxAssUni}, which corresponds to a graphical presentation of the $C \in AAC(\mathscr{O})$ in (\ref{AACproof}): The maximum number of objects in $C$ with sizes bigger than $2^{n-c(n)}$ (i.e. the number of \textbf{\textcolor{red}{red}} boxes above the level $c(n)$) is given by 
\begin{equation*}
    2^0 + 2^1 + \cdots + 2^{c(n)-1} = \sum_{j=0}^{c(n)-1}2^j = 2^{c(n)}-1.
\end{equation*}
On the other hand, the total number of the rest of the objects in $C$ with sizes smaller or equal than $2^{n-c(n)}$ (i.e. the number of \textbf{\textcolor{blue}{blue}} boxes) is bounded above by 
\begin{equation*}
\sum_{i=1}^{n-c(n)}\text{Card}(S(2^i)).
\end{equation*}
\begin{figure}[h]
    \centering
\includegraphics[width=0.5\textwidth]{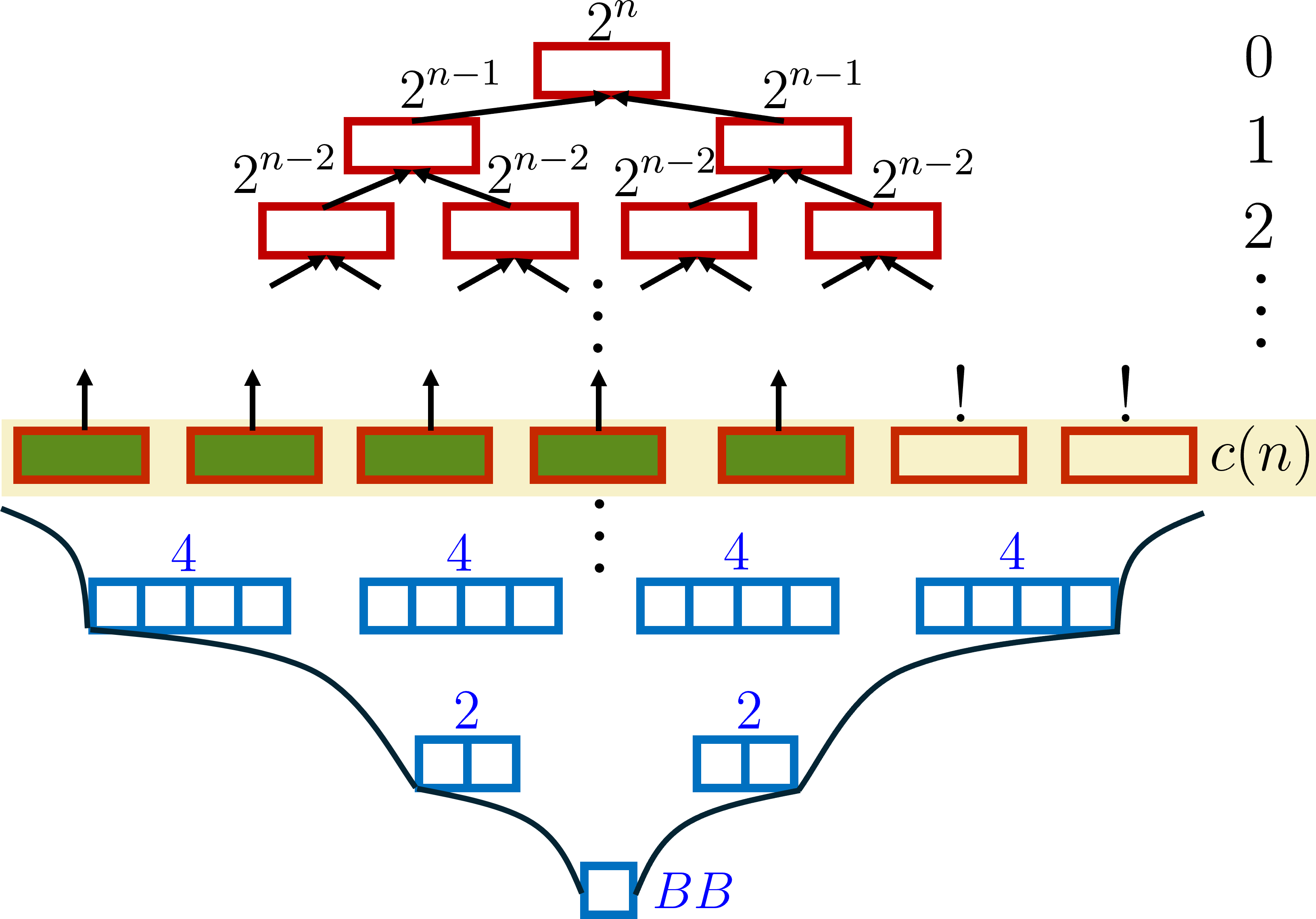}
    \caption{Graphical description of the construction of a maximal Assembly Path for a Binary Decomposable Object $\mathscr{O}$ of size $s=2^{n}$.}
    \label{fig:MaxAssUni}
\end{figure}
\end{remark}
The upper bound in Lemma \ref{lemmaBD} could be extended to any object $\mathscr{O} \in BD(S)$.  
\begin{theorem}\label{TheoremBD}
Let $(S,\circ,BB)$ be an Assembly Space and $\mathscr{O} \in BD(S)$. Then
    \begin{equation}\label{maxassuni}
    {a}(\mathscr{O}) \leq (H(s(\mathscr{O}))-1) + \sum_{i=1}^{n_1} \min \left\{ \sum_{j=1}^{H(s(\mathscr{O}))}  \lfloor 2^{n_j-i} \rfloor, \text{Card}(S(2^{i})) \right\}.
\end{equation} 
\end{theorem}
\begin{proof}
Similarly to the proof of Lemma \ref{lemmaBD}, the inequality (\ref{maxassuni}) follows by constructing $C \in AAC(\mathscr{O})$ such that $L(C)$ is less or equal than the RHS of (\ref{maxassuni}). Since $\mathscr{O} \in BD(S)$ there exist $\mathscr{P}_i$ for $i=1,\dots,H(s(\mathscr{O}))$ such that 
\begin{equation*}
        \mathscr{O} \in (\mathscr{P}_{H(s(\mathscr{O}))} \circ \cdots \circ  (\mathscr{P}_3 \circ (\mathscr{P}_{2} \circ \mathscr{P}_{1} \underbrace{)) \cdots)}_{H(s(\mathscr{O}))-1}
    \end{equation*}
with each $\mathscr{P}_i$ satisfying $s(\mathscr{P}_i)=2^{n_i}$ and decomposable as in (\ref{pieces}). According to the proof of Lemma \ref{lemmaBD}, for each $\mathscr{P}_i$ we could use the decomposition (\ref{pieces}) to construct a $C_i \in AAC(\mathscr{P}_i)$ as given by (\ref{AACproof}). Using these Assembly Addition Chains for each $\mathscr{P}_i$ we construct $C \in AAC(\mathscr{O})$ explicitly given by   
\begin{equation*}
 C \equiv   (C_1,C_2,\dots,C_{H(s(\mathscr{O}))},(\mathscr{P}_2 \circ \mathscr{P}_1),(\mathscr{P}_3 \circ (\mathscr{P}_2 \circ \mathscr{P}_1)),\dots,(\mathscr{P}_{H(s(\mathscr{O}))} \circ \cdots \circ  (\mathscr{P}_3 \circ (\mathscr{P}_{2} \circ \mathscr{P}_{1} \underbrace{)) \cdots)}_{H(s(\mathscr{O}))-1} ).
\end{equation*}
Since the number of objects of size $2^i$ in $C_j$ is given by $\lfloor 2^{n_j-i} \rfloor$ and the maximum number of elements of size $k$ in $S$ is given by $\text{Card}(S(k))$, there exists $C' \in AAC(\mathscr{O})$ (obtained from $C$ by removing possible repeated objects) such that  
\begin{equation*}
    a(\mathscr{O}) \leq L(C') \leq (H(s(\mathscr{O}))-1) + \sum_{i=1}^{n_1} \min \left\{ \sum_{j=1}^{H(s(\mathscr{O}))}  \lfloor 2^{n_j-i} \rfloor, \text{Card}(S(2^{i})) \right\},
\end{equation*}
concluding with the proof.
\end{proof}

\begin{remark}
    The RHS of inequality (\ref{maxassuni}) reduces to the RHS of inequality (\ref{ineqBDObj}) when $s(\mathscr{O})=2^n$, since in this case $H(s(\mathscr{O}))=1$ with $n_1=n$.     
\end{remark}

\begin{remark}
For Assembly Spaces $S$ such that $S=BD(S)$ and $\text{Card}(S(k))=1$ for $k \in \mathbb{Z}^+$ (i.e. there exists only one object of each size) e.g. $1$-Strings then (\ref{maxassuni}) reduces to (\ref{ACineq}).    
\end{remark}

For the cases where the cardinality of the level sets $S(k)$ increases monotonically with $k$, like for $j$-Strings, CCG and CPoly, the second term on the RHS of (\ref{maxassuni}) 
\begin{equation*}
\sum_{i=1}^{n_1-c(n_1)} \min \left\{ \sum_{j=1}^{H(s(\mathscr{O}))} \lfloor 2^{n_j-i} \rfloor, \text{Card}(S(2^i))  \right\} + \sum_{i=n_1-(c(n_1)-1)}^{n_1} \min \left\{ \sum_{j=1}^{H(s(\mathscr{O}))} \lfloor 2^{n_j-i} \rfloor, \text{Card}(S(2^i))  \right\}
\end{equation*}
reduces to 
\begin{equation*}
   \sum_{i=1}^{n_1-c(n_1)} \text{Card}(S(2^i))  + \sum_{i=n_1-(c(n_1)-1)}^{n_1} \min \left\{ \sum_{j=1}^{H(s(\mathscr{O}))} \lfloor 2^{n_j-i} \rfloor, \text{Card}(S(2^i))  \right\}, 
\end{equation*}
since by definition $2^{k} \geq \text{Card}(S(2^{n_1-k}))$ for $k=c(n_1),\dots,n_1-1,$ and consequently for these type of Assembly Spaces the inequality (\ref{maxassuni}) takes the form  
    \begin{equation}\label{maxassunialt}
\begin{split}
    {a}(\mathscr{O}) & \leq (H(s(\mathscr{O}))-1) + \sum_{i=1}^{n_1-c(n_1)} \text{Card}(S(2^i)) \\
    & \qquad + \sum_{k=1}^{c(n_1)} \min \left\{ \sum_{j=1}^{H(s(\mathscr{O}))} 2^{n_j-(n_1-c(n_1)+k)}, \text{Card}(S(2^{n_1-c(n_1)+k})) \right\}.
    \end{split}
\end{equation}

\begin{remark}
The intuitive idea behind the formula (\ref{maxassunialt}) is presented in Figure \ref{fig:MaxAssgen}: Following the construction presented in Remark \ref{remrk6} the objects of sizes $2^{n_1},2^{n_2},\dots,2^{n_H}$ are constructed. Since $c(n_1)>c(n_2)>\cdots>c(n_H)$ all the objects considered (from bottom to top) in the assembling sequence of the object of size $2^{n_1}$ up to the level $c(n_1)$ include all the objects considered (from bottom to top) in the assembling sequences of the objects of sizes $2^{n_k}$ up to the level $c(n_k)$ for $k=2,\dots,H$. Hence, the second term in (\ref{maxassunialt}) counts the number of steps involved in the longest Assembly Addition Chain of the objects of sizes $2^{n_1},\dots,2^{n_H}$ (from bottom to top) up to the levels $c(n_1),\dots,c(n_H)$, respectively. \\
\begin{figure}[h]
    \centering
\includegraphics[width=\textwidth]{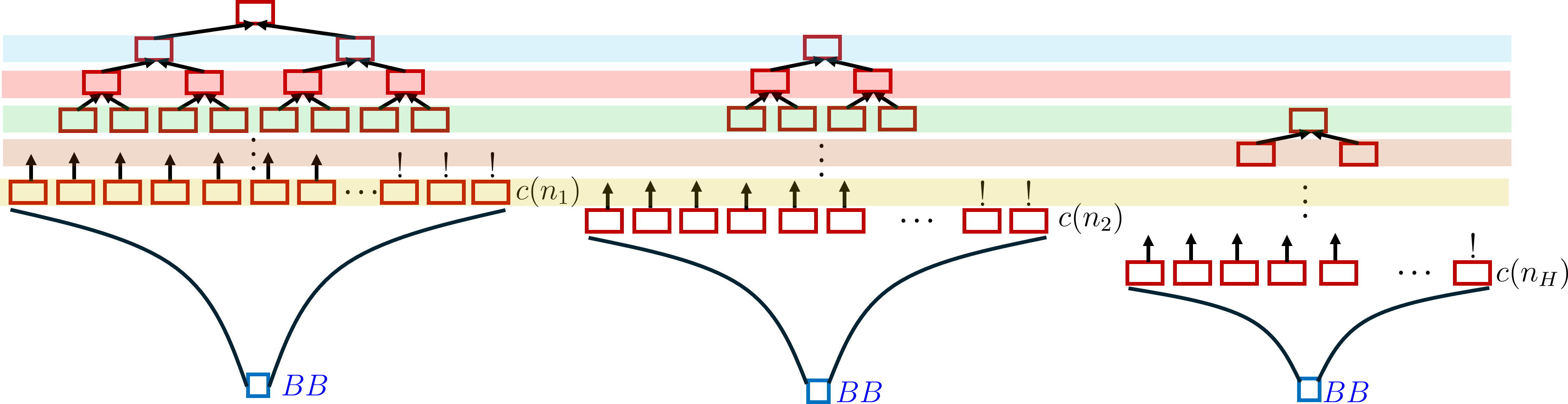}
    \caption{Graphical description of the construction of a Maximal Assembly Path for an object $\mathscr{O}$ of size $s=2^{n_1}+\cdots+2^{n_H}$ such that $n_1>\cdots>n_H$.}
    \label{fig:MaxAssgen}
\end{figure}

\noindent
The last term in (\ref{maxassunialt}) counts all the remaining binary steps (i.e. above the $c(n_k)$ level for the object of size $2^{n_k}$) of all the other maximal Assembly Paths avoiding over counting (since there cannot be more than  $\text{Card}(S(2^{k}))$ objects of size $2^{k}$) and the first term simply counts the extra $n_H(s)-1$ assembling steps required to glue together all the $n_H(s)$ objects.   
\end{remark}
\begin{remark}
    As expected, the RHS of inequality (\ref{maxassunialt}) reduces to the RHS of inequality (\ref{v2ineqBD}) when $s(\mathscr{O})=2^n$, since in this case $H(s(\mathscr{O}))=1$ with $n_1=n$.
\end{remark}

We denote by $G(s(\mathscr{O}))$ the RHS of inequality (\ref{maxassuni}) and define the function 
\begin{equation}\label{funBD}
\begin{split}
    Ma_{BD}:S & \to \mathbb{Z}^+ \\
    \mathscr{O} & \mapsto \min \{s(\mathscr{O})-1,G(s(\mathscr{O}))\}.
\end{split}    
\end{equation}
Then we have the following result for the Assembly Index of the Binary Decomposable objects in any Assembly Space.
\begin{corollary}
       Let $(S,\circ,BB)$ be an Assembly Space. For any $\mathscr{O} \in BD(S)$ the following inequalities hold
    \begin{equation}\label{UBBD}
    a(\mathscr{O}) \leq Ma_{BD}(\mathscr{O}) \leq s(\mathscr{O})-1.
    \end{equation}
\end{corollary}

In particular for Assembly Spaces for which $S=BD(S)$, e.g. $S=j$-Strings for any $j \in \mathbb{Z}^+$, (\ref{funBD}) provides a better upper bound than (\ref{ineq1}) for the Assembly Index of any object $\mathscr{O} \in S$. In the next section we will show explicitly how the upper bound given by $Ma_{BD}$ indeed follows the expected behavior from the maximum values display in Figure \ref{fig:allstrings}.





As indicated above there exist Assembly Spaces for which $BD(S) \subsetneq S$, and so inequality (\ref{UBBD}) does not hold for a generic element $\mathscr{O}$ of an arbitrary Assembly Space $S$. In particular, for $S=\text{CCG}$ (see Example \ref{nonexCSCPG}) and $S=\text{CPoly}$ (see Example \ref{nonexPoly}) the inclusion of $BD(S)$ in $S$ is strict. Nevertheless, the objects in these Assembly Spaces exhibit another type of decomposition that allows to derive also an upper bound better than $s(\mathscr{O})-1$. Moreover, they will also follow the behavior of Figure \ref{fig:allgraphs}.       
\begin{definition}\label{Def2PD}
    \textbf{(2-Piece decomposable object).} Let $(S,\circ,BB)$ be an Assembly Space. An object $\mathscr{O} \in S$ is called \textit{2-Piece decomposable} if both of the following conditions are satisfied:
    \begin{enumerate}[(i)]
        \item If there exists $\mathscr{B}_i \in BB$ with $i=1,\cdots,2k$ such that 
        \begin{equation}
            \mathscr{O} \in (\mathscr{B}_1 \circ (\mathscr{B}_2 \circ ( \cdots \circ (\mathscr{B}_{2k-1} \circ \mathscr{B}_{2k} \underbrace{) \cdots)  ))}_{2k-1},
        \end{equation}
        then there exists $\sigma \in S_{2k}$ such that
        \begin{equation}\label{2piece2k}
        \mathscr{O} \in ((\mathscr{B}_{\sigma(1)} \circ \mathscr{B}_{\sigma(2)}) \circ \cdots \circ ((\mathscr{B}_{\sigma(2k-3)} \circ \mathscr{B}_{\sigma(2k-2)}) \circ (\mathscr{B}_{\sigma(2k-1)} \circ \mathscr{B}_{\sigma(2k)}\underbrace{)) \cdots )}_{k}.
        \end{equation}
         \item If there exists $\mathscr{B}_i \in BB$ with $i=1,\cdots,2k+1$ such that
        \begin{equation}
 \mathscr{O} \in (\mathscr{B}_1 \circ (\mathscr{B}_2 \circ ( \cdots \circ (\mathscr{B}_{2k} \circ \mathscr{B}_{2k+1} \underbrace{) \cdots)  ))}_{2k},
        \end{equation}
        then there exists $\sigma \in S_{2k+1}$ such that
        \begin{equation}\label{2piece2k1}
            \mathscr{O} \in ((\mathscr{B}_{\sigma(1)} \circ \mathscr{B}_{\sigma(2)}) \circ \cdots \circ ((\mathscr{B}_{\sigma(2k-3)} \circ \mathscr{B}_{\sigma(2k-2)}) \circ (\mathscr{B}_{\sigma(2k-1)} \circ \mathscr{B}_{\sigma(2k)}\underbrace{)) \cdots )}_{k} \circ \mathscr{B}_{\sigma(2k+1)}.   
        \end{equation}
    \end{enumerate}   
\end{definition}
We denote by $2PD(S)$ the set of 2-piece decomposable objects of an Assembly space $S$. 

Defining the function 
\begin{equation}\label{func2pd}
\begin{split}
    Ma_{2PD}:S & \to \mathbb{Z}^+ \\
    \mathscr{O} & \mapsto \min \left\{ \bigg\lfloor \frac{s(\mathscr{O})}{2} \bigg\rfloor, \text{Card}(S(2))  \right\} + \bigg\lceil \frac{s(\mathscr{O})}{2} \bigg\rceil-1.
\end{split}    
\end{equation}
Then we have the following result for the Assembly Index of the 2-Piece decomposable objects in any Assembly Space.
\begin{theorem}
Let $(S,\circ,BB)$ be an Assembly Space. For any $\mathscr{O} \in 2PD(S)$ we have
\begin{equation}\label{ineq2PD}
    a(\mathscr{O}) \leq Ma_{2PD}(\mathscr{O}) \leq s(\mathscr{O})-1.
\end{equation}
\end{theorem}
\begin{proof}
Let $\mathscr{O} \in 2PD(S)$. If $s(\mathscr{O})=2k$ for $k \in \mathbb{Z}^+$, then by definition there exists $C \in AAC(\mathscr{O},BB)$ given in terms of (\ref{2piece2k}) by  
\begin{equation*}
\begin{split}
   C= & (\mathscr{B}_{\sigma(1)} \circ \mathscr{B}_{\sigma(2)}, \dots,\mathscr{B}_{\sigma(2k-1)} \circ \mathscr{B}_{\sigma(2k)}, (\mathscr{B}_{\sigma(2k-3)} \circ \mathscr{B}_{\sigma(2k-2)}) \circ (\mathscr{B}_{\sigma(2k-1)} \circ \mathscr{B}_{\sigma(2k)}),\dots,\\& ((\mathscr{B}_{\sigma(1)} \circ \mathscr{B}_{\sigma(2)}) \circ \cdots \circ ((\mathscr{B}_{\sigma(2k-3)} \circ \mathscr{B}_{\sigma(2k-2)}) \circ (\mathscr{B}_{\sigma(2k-1)} \circ \mathscr{B}_{\sigma(2k)}\underbrace{)) \cdots )}_{k}). 
   \end{split}
\end{equation*}
Since the first $k$ objects in the sequence have size two and the number of objects of size two in $S$ is $\text{Card}(S(2))$, then there exists $C' \in AAC(\mathscr{O},BB)$ (obtained from $C$ by removing possible repeated objects of size two) such that  
\begin{equation*}
    L(C')= \min \{k, \text{Card}(S(2)) \} + k-1.
\end{equation*}
Analogously, if $s(\mathscr{O})=2k+1$ for some $k \in \mathbb{Z}^+$, then there exists $C \in AAC(\mathscr{O},BB)$ given implicitly by (\ref{2piece2k1})
\begin{equation*}
\begin{split}
   C= & (\mathscr{B}_{\sigma(1)} \circ \mathscr{B}_{\sigma(2)}, \dots,\mathscr{B}_{\sigma(2k-1)} \circ \mathscr{B}_{\sigma(2k)}, (\mathscr{B}_{\sigma(2k-3)} \circ \mathscr{B}_{\sigma(2k-2)}) \circ (\mathscr{B}_{\sigma(2k-1)} \circ \mathscr{B}_{\sigma(2k)}),\dots,\\& ((\mathscr{B}_{\sigma(1)} \circ \mathscr{B}_{\sigma(2)}) \circ \cdots \circ ((\mathscr{B}_{\sigma(2k-3)} \circ \mathscr{B}_{\sigma(2k-2)}) \circ (\mathscr{B}_{\sigma(2k-1)} \circ \mathscr{B}_{\sigma(2k)}\underbrace{)) \cdots )}_{k},\\
   & ((\mathscr{B}_{\sigma(1)} \circ \mathscr{B}_{\sigma(2)}) \circ \cdots \circ ((\mathscr{B}_{\sigma(2k-3)} \circ \mathscr{B}_{\sigma(2k-2)}) \circ (\mathscr{B}_{\sigma(2k-1)} \circ \mathscr{B}_{\sigma(2k)}\underbrace{)) \cdots )}_{k}\circ \mathscr{B}_{2k+1}) 
   \end{split}
\end{equation*}
So, there exists $C' \in AAC(\mathscr{O},BB)$ (obtained from $C$ by removing possible repeated objects of size two) such that 
\begin{equation*}
    L(C')= \min \{k, \text{Card}(S(2)) \} + k.
\end{equation*}
The last inequality in (\ref{ineq2PD}) simply follows from 
\begin{equation*}
   \bigg\lfloor \frac{s(\mathscr{O}}{2} \bigg\rfloor + \bigg\lceil \frac{s(\mathscr{O})}{2} \bigg\rceil = s(\mathscr{O}).
\end{equation*}
\end{proof}


In particular for Assembly Spaces for which $2PD(S)=S$,(\ref{func2pd}) provides a better upper bound than (\ref{ineq1}) for the Assembly Index of any object $\mathscr{O} \in S$. In the next section we will show explicitly how these upper bounds could be computed for $S=\text{CCG}$ and $S=\text{CPoly}$.




\section{Examples}\label{secExa}
In the previous section we introduced the functions $Ma_{BD}$ (\ref{funBD}) and $Ma_{2PD}$ (\ref{func2pd}) and proved that these provide upper bounds, better than $s(\mathscr{O})-1$, for the Assembly Index of Binary Decomposable and 2-Piece Decomposable objects, respectively. In this section, we show explicitly how $Ma_{BD}$ gives an upper bound for the AI of any string in $S=j$-String (for any $j \in \mathbb{Z}^+$) and how $Ma_{2PD}$ gives an upper bound for the AI of any graph in $S=\text{CCG}$ and any Polyomino in $S=\text{CPoly}$.             

\subsection{Addition Chains of $j$-Strings}
Since any $j$-String of size $L$ could be naturally splitted into sub-strings of size $k$ and $L-k$ for any $1 \leq k<L$ and the Multi-Magma structure considered in this paper is concatenation (see Example \ref{exmpstrings}), it follows that for this Assembly Space, $BD(S)=S$.  Hence, as explained in the previous section (\ref{funBD}) gives an upper bound for the AI of any $j$-String for any $j \in \mathbb{Z}^+$. Since the cardinalities of the sets $S(2^i)$ are different for directed and undirected strings, we consider each case separately.      

\subsubsection{Directed Strings}
Let $j \in \mathbb{Z}^+$. Since 
\begin{equation}\label{carddir}
\text{Card}(j\text{-String}_d(2^i))=j^{2^i}  
\end{equation}
for any $i \in \mathbb{Z}_{\geq 0}$, then according to (\ref{UBBD}) an upper bound for the AI of a directed $j$-string $\mathscr{O}$ of size $2^{n_1}+2^{n_2}+\cdots + 2^{n_H}$ (with $n_1>n_2>\cdots > n_H$) is 
\begin{equation}\label{maxassstring}
\begin{split}
   {Ma}^{j\text{-Strings}_d}(\mathscr{O}) = \min \left\{ s-1,(H(s)-1) + \sum_{i=1}^{n_1-c(n_1)} j^{2^i}  + \sum_{k=1}^{c(n_1)} \min \left( \sum_{i=1}^{H} 2^{n_i-(n_1-c(n_1)+k)}, j^{(2^{(n_1-c(n_1)+k)})} \right) \right\}
    \end{split}
\end{equation}
with $c(n)$ in (\ref{c_n}) given in terms of the function
\begin{equation}
    \mu^{j\text{-Strings}_d}_n(i)=2^{i}-j^{2^{n-i}}.
\end{equation}

\begin{figure}[h]
    \centering
\includegraphics[width=0.6\textwidth]{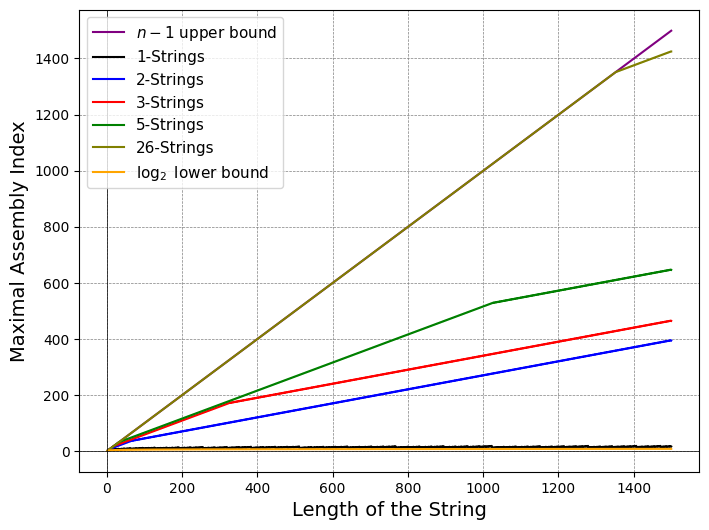}
    \caption{Maximum Assembly Index of $j$-Strings for $j=1,2,3,5,26$.}
    \label{fig:MaxAssString1}
\end{figure}

In Figure \ref{fig:MaxAssString1} the values of the Maximal Assembly Index function for directed $j$-Strings (\ref{maxassstring}) are presented for the cases $j=1,2,3,5$ and $26$. From this figure it is apparent that up to certain size (which depends monotonously on $j$) the MAI coincides with the upper bound $s(\mathscr{O})-1$ and then $Ma^{j\text{-Strings}_d}$ is \textit{strictly} less than this upper bound, exhibiting a \textit{quasi}-linear (saw-tooth shaped curve) monotonous increasing behavior, which is precisely the one presented by the maximum AI in Figure \ref{fig:allstrings}.

\subsubsection{Undirected}

\begin{figure}[h!]
    \centering 
\begin{minipage}[t]{.45\textwidth}
\begin{subfigure}{\textwidth}
  \includegraphics[width=\linewidth]{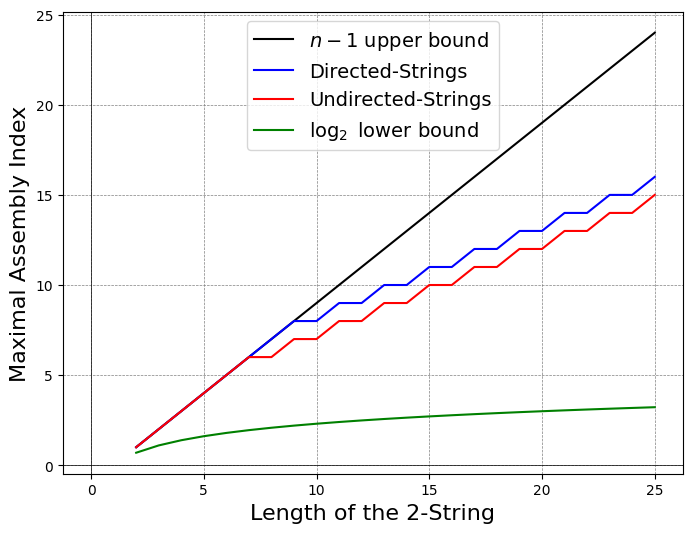}
  \caption{j=2}
  \label{fig:1}
\end{subfigure}\hfil 
\begin{subfigure}{\textwidth}
  \includegraphics[width=\linewidth]{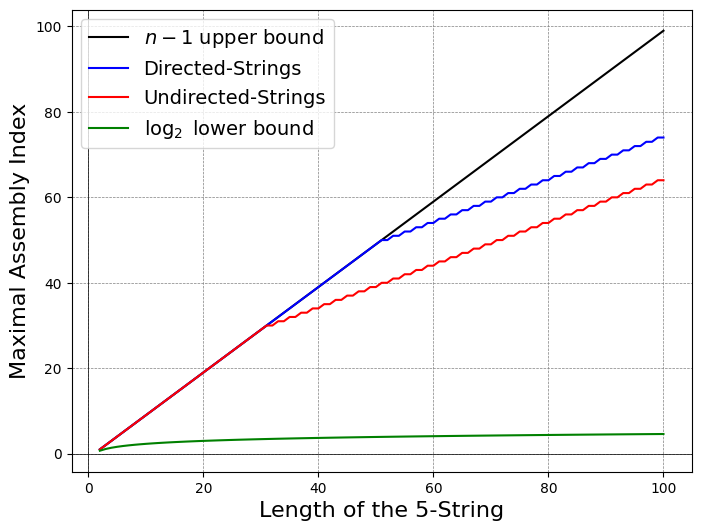}
  \caption{j=5}
  \label{fig:3}
\end{subfigure}
\end{minipage}\hfil
\begin{minipage}[t]{.45\textwidth}
\begin{subfigure}{\textwidth}
  \includegraphics[width=\linewidth]{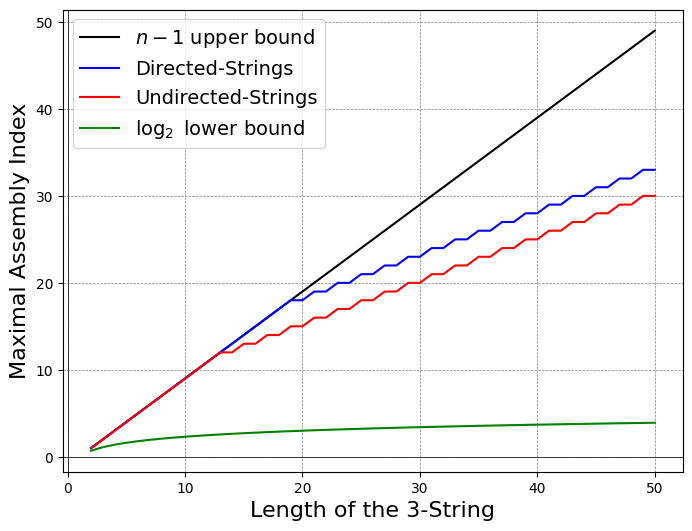}
  \caption{j=3}
  \label{fig:2}
\end{subfigure}\hfil 
\begin{subfigure}{\textwidth}
  \includegraphics[width=\linewidth]{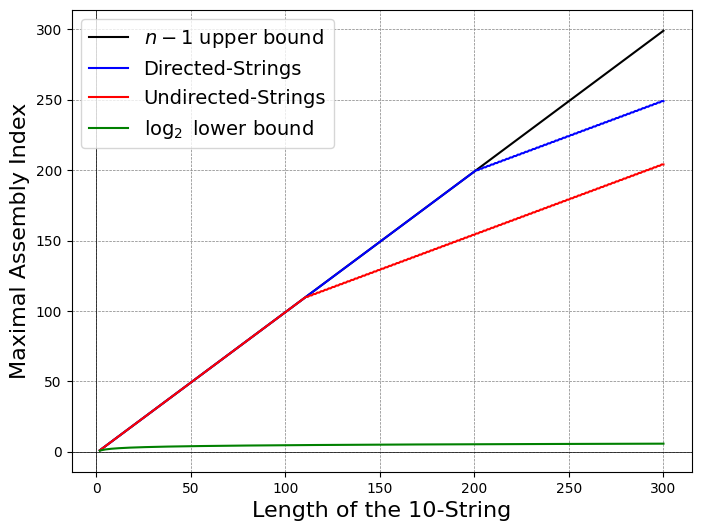}
  \caption{j=10}
  \label{fig:4}
\end{subfigure}\hfil 
\end{minipage}
\caption{Comparison of the MAI for undirected (\textbf{\textcolor{red}{red}}) and directed (\textbf{\textcolor{blue}{blue}}) $j$-Strings for $j=2,3,5,10$. Also the universal upper $\ell(\mathscr{O})-1$ (\textbf{\textcolor{black}{black}}) and lower $\floor{\log_2(s(\mathscr{O}))}$ (\textbf{\textcolor{green}{green}}) bounds for the AI are shown.}
\label{DandUStrings}
\end{figure}

On the other hand, since the number of undirected $j$-Strings of length $2^i$ is 
\begin{equation}\label{cardund}
         \text{Card}(j\text{-String}_u(2^i)) =\underbrace{\frac{j^{2^i}-j^{2^{i-1}}}{2}}_{\text{non-palindromes}}+\underbrace{j^{2^{i-1}}}_{\text{palindromes}} = \frac{j^{2^i}+j^{2^{i-1}}}{2}
\end{equation}
for any $i \in \mathbb{Z}_{\geq 0}$, in this case for an undirected $j$-string $\mathscr{O}$ of size $s=2^{n_1}+2^{n_2}+\cdots + 2^{n_H}$ (with $n_1>n_2>\cdots > n_H$) the function (\ref{funBD}) takes the form 
\begin{equation}\label{maxassustring}
\begin{split}
    Ma^{j\text{-Strings}_u}(\mathscr{O})  & = \min \bigg\{ s-1 ,(H(s)-1)  +\sum_{i=1}^{n_1-c(n_1)} \frac{j^{2^{n-i}}+j^{2^{n-i-1}}}{2} \\
    & \qquad + \sum_{k=1}^{c(n_1)} \min \left( \sum_{i=1}^{H} 2^{n_i-(n_1-c(n_1)+k)}, \frac{j^{(2^{(n_1-c(n_1)+k)})}+j^{(2^{(n_1-c(n_1)+k-1)})}}{2} \right) \bigg\},
\end{split}
\end{equation}
with $c(n)$ in (\ref{c_n}) now given by the function
\begin{equation}
    \mu^{j\text{-Strings}_u}_n(i)=2^{i}-\frac{j^{2^{n-i}}+j^{2^{n-i-1}}}{2} = \frac{2^{i+1}-(j^{2^{n-i}}+j^{2^{n-i-1}})}{2}.
\end{equation}
In Figure \ref{DandUStrings} the comparisons of the values of the upper bounds for undirected (\textbf{\textcolor{red}{red}}) (\ref{maxassustring}) and directed (\textbf{\textcolor{blue}{blue}}) $j$-Strings (\ref{maxassustring}) are presented for the cases $j=2,3,5$ and $10$. From this figure it is evident that the  behavior of both upper bounds is the same, but the one of the undirected Strings is always lower or equal than the one of the directed Strings, which follows from the fact that the cardinalities (\ref{carddir}) and (\ref{cardund}) satisfy the inequality
\begin{equation*}
   \text{Card}(j\text{-String}_d(2^k)) \geq  \text{Card}(j\text{-String}_u(2^k)),
\end{equation*}
for every $j,k \in \mathbb{Z}^+$.

\subsection{Addition Chains of Colored Connected Graphs}
\begin{figure}[h]
    \centering
\includegraphics[width=0.6\textwidth]{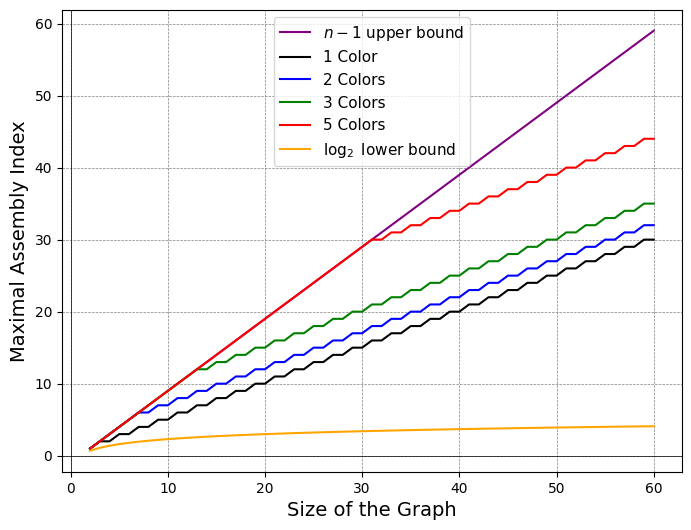}
    \caption{Maximum Assembly Index of CCG for $1,2,3$ and $5$ colors.}
    \label{fig:MaxAsSGraphs}
\end{figure}

For any (colored) Connected graph it is always possible (i) to remove one edge such that the resulting graph is still connected and (ii) if it has an even number of edges then it admits a decomposition into paths of length two (this is a well-known result by Kotzig, see e.g. \cite{bondy1990perfect} and the original paper \cite{kotzig1957theorie}). Hence, equipped with the Multi-Magma structure from Example \ref{exmgraphs}, the Assembly Space $S=\text{CCG}$ is such that $S=2PD(S)$. Consequently, function (\ref{func2pd}) provides an upper bound for the AI of all objects in $S=\text{CCG}$.  \\\\
For CCG the function \ref{func2pd} is given by
\begin{equation}
 Ma^{\text{CCG;}colors}(\mathscr{O}) = \min \left\{ \bigg\lfloor \frac{s(\mathscr{O})}{2} \bigg\rfloor, \binom{colors}{2}  \right\} + \bigg\lceil \frac{s(\mathscr{O})}{2} \bigg\rceil-1.   
\end{equation}
In Figure \ref{fig:MaxAsSGraphs} the values of the Maximal Assembly Index function for CSCPGs are presented for the cases $colors=1,2,3$ and $5$. From this figure, similar to the case of $j$-Strings, it is clear to realize that up to a certain size (which depends monotonously on $colors$) the MAI coincides with the upper bound $s(\mathscr{O})-1$ and then the MAI is \textit{strictly} less than this upper bound, exhibiting a \textit{quasi}-linear (saw-tooth shaped curve) monotonous increasing behavior. This is precisely the same behavior as in Figure \ref{fig:allgraphs}.

\subsection{Addition Chains of Colored Polyominoes}
\begin{figure}[h]
    \centering
\begin{tikzpicture}[scale=1, every node/.style={inner sep=1pt}]
\tikzstyle{rededge}=[line width=1pt, red]
\tikzstyle{blueedge}=[line width=1pt, lightblue]
\tikzstyle{blackedge}=[line width=1pt]


\node at (-1,1.5) {$\bigg\{\, 
\begin{tikzpicture}[scale=0.8]
  \coordinate (a) at (0,0);
    \filldraw[lightblue] (a) circle (1.5pt);
  
\end{tikzpicture}
\bigg\}\; \circ \; \bigg\{ \rule{0pt}{12pt}$};

\begin{scope}[shift={(0.1,3)}]
  \coordinate (a) at (0,0);
  \coordinate (b) at (0,0.5);
  \coordinate (c) at (0,-0.5);
  \draw[blackedge] (a)--(b)--(c);
    \foreach \pt in {b,c} {
    \filldraw[black] (\pt) circle (1.5pt);
  }
  \draw[fill=white, draw=black] (a) circle (1.8pt);
\end{scope}

\node at (0.9,0) {$\,\bigg\}  \;=\; \bigg\{ \rule{0pt}{12pt}$};

\begin{scope}[shift={(1.8,2.5)}]
  \coordinate (a) at (0,0);
  \coordinate (b) at (0,0.5);
  \coordinate (c) at (0,1.0);
  \draw[blackedge] (a)--(b)--(c);
    \foreach \pt in {a,c} {
    \filldraw[black] (\pt) circle (1.2pt);
  }
  \filldraw[red] (b) circle (1.2pt);
\end{scope}


\node at (2.1,3) {$\bigg\} \rule{0pt}{12pt}$};

\node at (-1,3) {$\bigg\{\, 
\begin{tikzpicture}[scale=0.8]
  \coordinate (a) at (0,0);
    \filldraw[red] (a) circle (1.5pt);
  
\end{tikzpicture}
\bigg\}\; \circ \; \bigg\{ \rule{0pt}{12pt}$};

\begin{scope}[shift={(0.1,1.5)}]
  \coordinate (a) at (0,0);
  \coordinate (b) at (0,0.5);
  \coordinate (c) at (0,-0.5);
  \draw[blackedge] (a)--(b)--(c);
    \foreach \pt in {a,c} {
    \filldraw[black] (\pt) circle (1.5pt);
  }
  \draw[fill=white, draw=black] (b) circle (1.8pt);
\end{scope}

\node at (0.9,1.5) {$\,\bigg\}  \;=\; \bigg\{ \rule{0pt}{12pt}$};

\begin{scope}[shift={(1.8,1)}]
  \coordinate (a) at (0,0);
  \coordinate (b) at (0,0.5);
  \coordinate (c) at (0,1.0);
  \draw[blackedge] (a)--(b)--(c);
    \foreach \pt in {a,b} {
    \filldraw[black] (\pt) circle (1.2pt);
  }
  \filldraw[lightblue] (c) circle (1.2pt);
\end{scope}


\node at (2.1,1.5) {$\bigg\} \rule{0pt}{12pt}$};

\node at (-1.7,0) {$\bigg\{\, \rule{0pt}{12pt}$};
\begin{scope}[shift={(-1.45,0)}]
  \coordinate (a) at (0,0);
  \coordinate (b) at (0,0.5);
  \coordinate (c) at (0,-0.5);
  \draw[blackedge] (a)--(b)--(c);

    \filldraw[black] (c) circle (1.5pt);
    \filldraw[red] (b) circle (1.5pt);
  \draw[fill=white, draw=black] (a) circle (1.8pt);
\end{scope}
\node at (-0.7,0) {$ \bigg\}\; \circ \; \bigg\{ \rule{0pt}{12pt}$};

\begin{scope}[shift={(0.1,0)}]
  \coordinate (a) at (0,0);
  \coordinate (b) at (0,0.5);
  \coordinate (c) at (0,-0.5);
  \draw[blackedge] (a)--(b)--(c);
    \filldraw[green] (c) circle (1.5pt);
    \filldraw[lightblue] (a) circle (1.5pt);
  \draw[fill=white, draw=black] (b) circle (1.8pt);
\end{scope}

\node at (0.9,3) {$\,\bigg\}  \;=\; \bigg\{ \rule{0pt}{12pt}$};

\begin{scope}[shift={(1.8,-0.5)}]
  \coordinate (a) at (0,0);
  \coordinate (b) at (0,0.5);
  \coordinate (c) at (0,1.0);
  \coordinate (d) at (0.5,1.0);
  \coordinate (e) at (1.0,1.0);
  \draw[blackedge] (a)--(b)--(c);
  \draw[blackedge] (c)--(d)--(e);
  \draw[fill=white, draw=black] (b) circle (1.8pt);
    \filldraw[black] (a) circle (1.2pt);
  \filldraw[red] (c) circle (1.2pt);
  \filldraw[lightblue] (d) circle (1.2pt);
  \filldraw[green] (e) circle (1.2pt);
\end{scope}
\node at (3.1,0) {\Huge$ , \rule{0pt}{12pt}$};

\begin{scope}[shift={(3.4,-0.5)}]
  \coordinate (a) at (0,0);
  \coordinate (b) at (0,0.5);
  \coordinate (c) at (0,1.0);
  \coordinate (d) at (0.5,0);
  \coordinate (e) at (1.0,0);
  \draw[blackedge] (a)--(b)--(c);
  \draw[blackedge] (a)--(d)--(e);
  \draw[fill=white, draw=black] (b) circle (1.8pt);
    \filldraw[black] (a) circle (1.2pt);
  \filldraw[red] (c) circle (1.2pt);
  \filldraw[lightblue] (d) circle (1.2pt);
  \filldraw[green] (e) circle (1.2pt);
\end{scope}
\node at (4.7,0) {\Huge$ , \rule{0pt}{12pt}$};

\begin{scope}[shift={(5.5,-0.5)}]
  \coordinate (a) at (0,0);
  \coordinate (b) at (0,0.5);
  \coordinate (c) at (0,1.0);
  \coordinate (d) at (0.5,0.5);
  \coordinate (e) at (-0.5,0.5);
  \draw[blackedge] (a)--(b)--(c);
  \draw[blackedge] (d)--(b)--(e);
  \draw[fill=white, draw=black] (d) circle (1.8pt);
    \filldraw[black] (a) circle (1.2pt);
  \filldraw[red] (c) circle (1.2pt);
  \filldraw[lightblue] (b) circle (1.2pt);
  \filldraw[green] (e) circle (1.2pt);
\end{scope}
\node at (6.4,0) {\Huge$ , \rule{0pt}{12pt}$};

\begin{scope}[shift={(6.7,-0.5)}]
  \coordinate (a) at (0,0);
  \coordinate (b) at (0,0.5);
  \coordinate (c) at (0,1.0);
  \coordinate (d) at (0.5,0.5);
  \coordinate (e) at (1.0,0.5);
  \draw[blackedge] (a)--(b)--(c);
  \draw[blackedge] (d)--(b)--(e);
  \draw[fill=white, draw=black] (e) circle (1.8pt);
    \filldraw[black] (a) circle (1.2pt);
  \filldraw[red] (c) circle (1.2pt);
  \filldraw[lightblue] (d) circle (1.2pt);
  \filldraw[green] (b) circle (1.2pt);
\end{scope}

\node at (8,0) {$\bigg\} \rule{0pt}{12pt}$};

\end{tikzpicture}
    \caption{Multi-Magma structure over the set of Skeletons of Polyominoes.}
    \label{fig:MultiSke}
\end{figure}
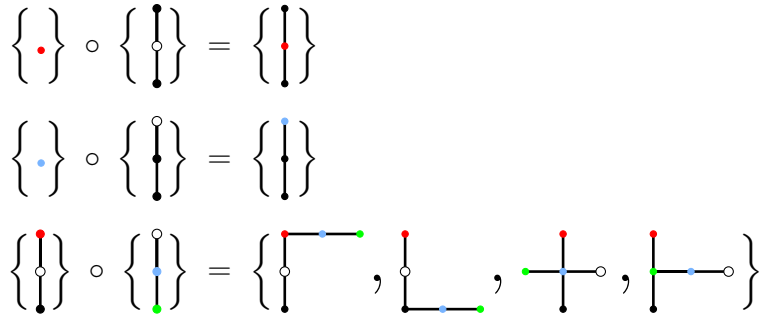

For any colored Polyomino we could attach an unique vertex-colored planar graph, such that a colored square corresponds to a colored vertex and a shared edge between two squares corresponds to an edge connecting the corresponding vertices. This well-defined mapping from the set of colored Polyominoes to the set of vertex-colored Connected Graphs is bijective and provides a representation of any colored Polyomino of size $n$ (using $m$ colors) by a vertex-colored Connected Graph with $n-1$ or $n$ edges (using $m$ colors), which we called the \textit{Skeleton} of the Polyomino. \\\\
Since the construction of the skeleton amounts to the construction of the Polyomino, we define a Multi-Magma structure over the set of colored Polyominoes (different from Example \ref{CpolyMM}) in terms of the Skeletons: The Building Blocks correspond to (a) colored vertices and (b) hooked-graphs with two edges and colored vertices. See Figure \ref{fig:MultiSke} for the Multi-Magma structure over the Set of Skeletons for colored Polyominoes. 

\begin{figure}[h]
\centering
\begin{subfigure}[b]{0.6\textwidth}
\includegraphics[width=\linewidth]{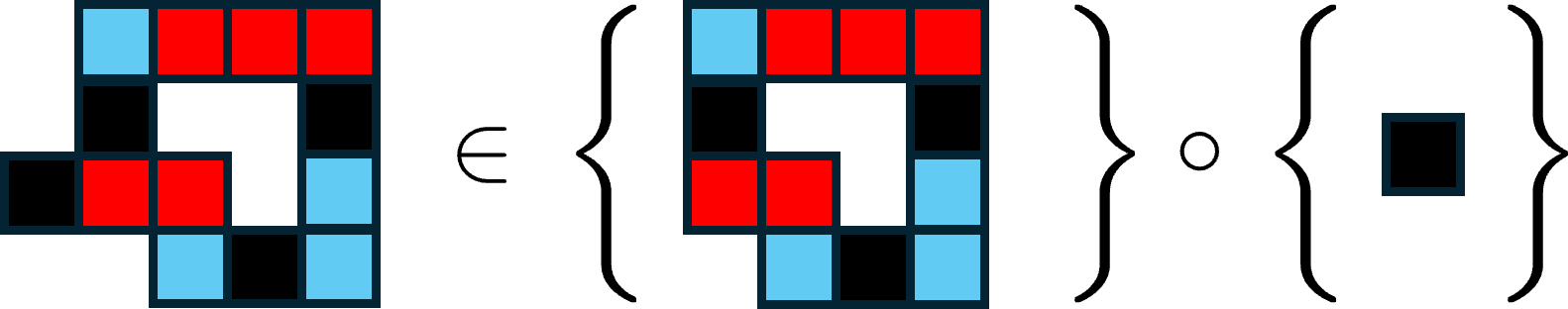}
\caption{For any polyomino of odd number of blocks we can remove one block such that the resulting polyomino is connected and contains an even number of blocks.\medskip\medskip}
\label{subfig:image1}
\end{subfigure}\hfill
\begin{subfigure}[b]{0.3\textwidth}
\includegraphics[width=\linewidth]{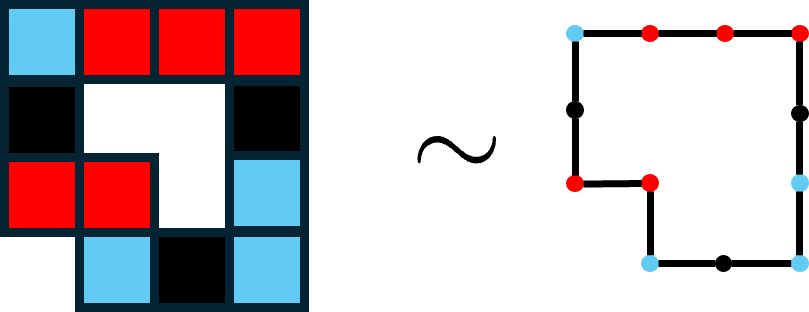}
\caption{Any polyomino has a skeleton graph representation that can be built using colored vertices and hooked-graphs.}
\label{subfig:image2}
\end{subfigure}\hfill
\begin{subfigure}[b]{0.65\textwidth}
\includegraphics[width=\linewidth]{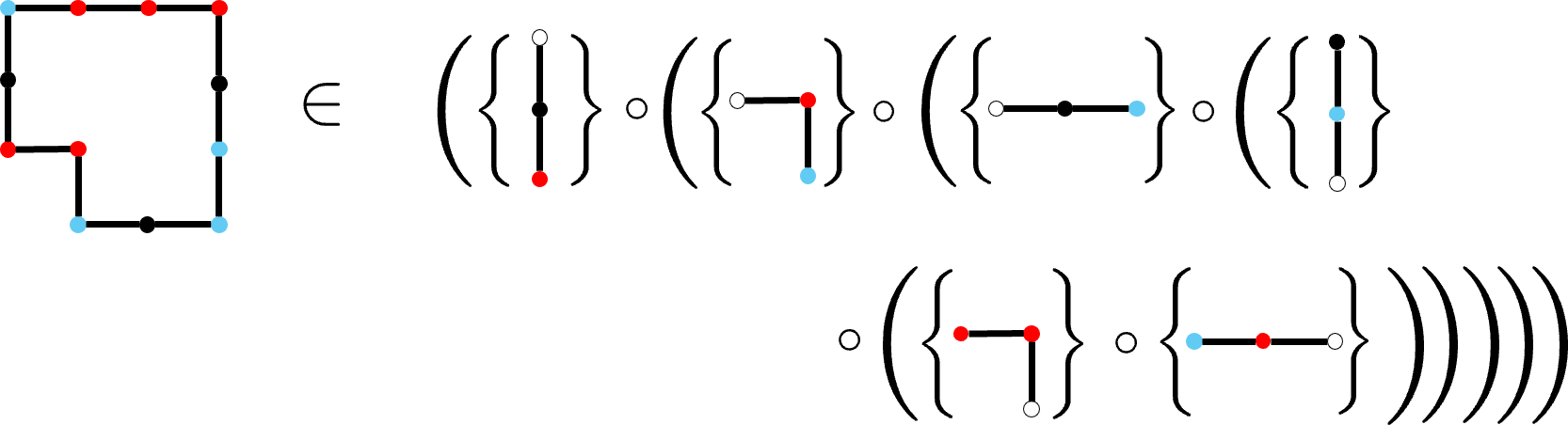}
\caption{Two-Piece decomposition of the skeleton graph representation of a polyomino.}
\label{subfig:image3}
\end{subfigure}
\caption{Procedure to obtain a two-piece decomposition of any colored polyomino. a) Conversion to even number of blocks, b) skeleton representation, c) 2-Piece decomposition of skeleton.}
\label{fig:exskel}
\end{figure}

As an example of how this Multi-Magma structure could be used to produce Assembly Addition Chains for colored Polyominoes, in Figure \ref{fig:exskel} we present as an example an Assembly Addition Chain for the polyomino in Figure \ref{figbdpoly} (a).  

\begin{figure}[!h]
    \centering
\includegraphics[width=0.6\textwidth]{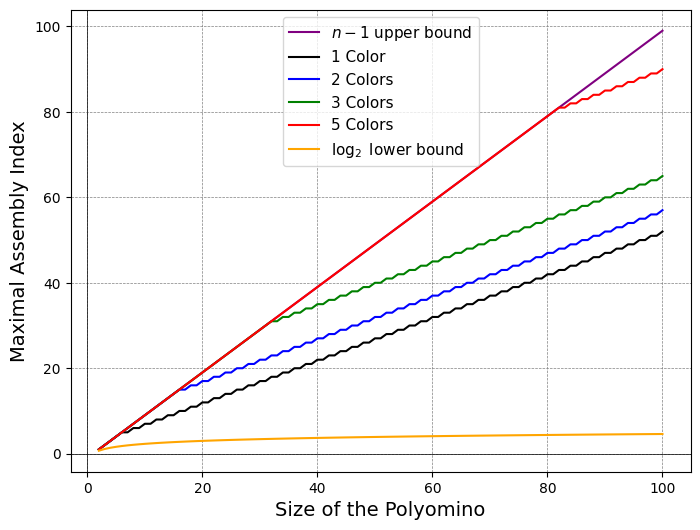}
    \caption{Maximum Assembly Index of Polyominoes for $1,2,3,4$ colors.}
    \label{fig:MaxAsSPoly}
\end{figure}

Using this Multi-Magma structure and the length two path-decomposition of Connected Graphs, we derive that the Assembly Index of any $\mathscr{O} \in \text{CPoly}$ is bounded above by
\begin{equation}
    Ma^{\text{Poly;}colors}(\mathscr{O}) = \min \left\{ \bigg\lfloor \frac{s(\mathscr{O})-1}{2} \bigg\rfloor, 2 \left(colors+\binom{colors}{2} \right)  \right\} + \bigg\lceil \frac{s(\mathscr{O})-1}{2} \bigg\rceil -1 + \underbrace{1}_{\text{First color fixed}},
\end{equation}
where $s(\mathscr{O})-1$ is the number of edges of the skeleton of the Polyomino $\mathscr{O}$, the second argument in the $\min$ function follows from the fact that there are $color+ \binom{color}{2}$ colorings for the hooked-graphs with two edges of type 
\begin{tikzpicture}[scale=0.5,baseline=7.5ex]
\tikzstyle{rededge}=[line width=1pt, red]
\tikzstyle{blueedge}=[line width=1pt, lightblue]
\tikzstyle{blackedge}=[line width=1pt]
\begin{scope}[shift={(0.1,3)}]
  \coordinate (a) at (0,0);
  \coordinate (b) at (0,0.5);
  \coordinate (c) at (0,-0.5);
  \draw[blackedge] (a)--(b)--(c);
    \foreach \pt in {b,c} {
    \filldraw[black] (\pt) circle (1.8pt);
  }
  \draw[fill=white, draw=black] (a) circle (1.8pt);
\end{scope}
\end{tikzpicture}
and $color + color(color-1)$ colorings for the hooked-graphs with two edges of type 
\begin{tikzpicture}[scale=0.5, baseline=-6.0ex]
\tikzstyle{rededge}=[line width=1pt, red]
\tikzstyle{blueedge}=[line width=1pt, lightblue]
\tikzstyle{blackedge}=[line width=1pt]
\begin{scope}[shift={(2,-2)}]
  \coordinate (a) at (0,0);
  \coordinate (b) at (0,0.5);
  \coordinate (c) at (0,-0.5);
  \draw[blackedge] (a)--(b)--(c);
    \foreach \pt in {a,c} {
    \filldraw[black] (\pt) circle (1.8pt);
  }
  \draw[fill=white, draw=black] (b) circle (1.8pt);
\end{scope}
\end{tikzpicture} and the +1 corresponds to fixing the color of the square associated to the colorless vertex of the first hooked-graph with two edges used in the construction of the skeleton (i.e. one of the first two gluing operations displayed in Figure \ref{fig:MultiSke}). \\\\
In Figure \ref{fig:MaxAsSPoly} the values of the Maximal Assembly Index function for $S=\text{Poly}$ are presented for the cases $\text{colors}=1,2,3$ and $4$. From this figure, similar to the case of $j$-Strings and CSCG, it is clear to realize that up to a certain size (which depends monotonously on $colors$) the MAI coincides with the upper bound $s(\mathscr{O})-1$ and then the MAI is \textit{strictly} less than this upper bound, exhibiting a \textit{quasi}-linear (saw-tooth shaped curve) monotonous increasing behavior.

\section{Outlook and Conclusions}\label{secOut}

In this paper we introduce the concept of Assembly Addition Chains over an arbitrary set $S$, extending the standard  notion of Addition Chains over $\mathbb{Z}^+$. The central algebraic structure for this extension is the Multi-Magma, which formalizes the intuitive notion of gluing/assembling objects of a set, such that the standard semi-group structure $(\mathbb{Z}^+,+,1)$ is replaced by the Assembly Multi-Magma structure $(S,\circ,BB)$. When (i) every object $\mathscr{O} \in S$ is guaranteed to be constructed via $\circ$ from  the subset $BB$ and (ii) $\circ$ is compatible with the size function (\ref{sizeprop}), then $(S,\circ,BB)$ becomes an Assembly Space. \\\\
Analogously to the way $AC(n)$ is defined over $(\mathbb{Z}^+,+,1)$ for any $n \in \mathbb{Z}^+$, the set of Assembly Addition Chains $AAC(\mathscr{O},BB)$ of the object $\mathscr{O} \in S$  over any Assembly Space $(S,\circ,BB)$ is defined as the set of sequences $C$ of objects in $S$ such that (i) the last object in $C$ is $\mathscr{O}$ and (ii) every object in the sequence is the gluing of (a) two predecessor objects in the sequence, or (b) two Building Blocks, or (c) one Building Block and one predecessor object in the sequence. Moreover, the set of Optimal Assembly Addition Chains of the object $\mathscr{O} \in S$, $OAAC(\mathscr{O},BB)$, corresponds to the set of the shortest Assembly Addition Chains of $\mathscr{O}$ and their length is the so-called Assembly Index $a(\mathscr{O})$ of the object $\mathscr{O} \in S$.      \\\\
For any $n \in \mathbb{Z}^+$, the set $AC(n)$ and the number $\ell(n)$ correspond to the set of shortest additive presentations of $n$ and to the additive integer complexity of $n$, respectively. In the same fashion for any object $\mathscr{O} \in S$ in an Assembly Space $(S,\circ,BB)$, the set $OAAC(\mathscr{O})$ corresponds to the set of shortest assembly paths/presentations of $\mathscr{O}$ and the number $a(\mathscr{O})$ is defined as the assembly complexity of the object. \\\\
Assembly Addition Chains extend the well-studied field of standard Addition Chains, providing a formal way to represent any object $\mathscr{O}$ as a sequence $C \in AAC(\mathscr{O},BB)$ of recursive gluing operations and the Assembly Index $a(\mathscr{O})$ defines the complexity involved in the optimal recursive representation(s) of $\mathscr{O}$. As the representations given by $OAC(n)$ result fundamental for optimal computer algebra, it is expected that the elements in $OAAC(\mathscr{O})$ become fundamental regarding the compressibility and optimal computer representations of the objects $\mathscr{O} \in S$.   \\\\
The computation of $a(\mathscr{O})$ for an arbitrary object $\mathscr{O} \in S$ is an NP-complete problem and consequently the derivation of bounds is fundamentally advantageous for estimations. In this paper we proved that for any $\mathscr{O} \in S=j$-Strings, $\ell(s(\mathscr{O})) \leq a(\mathscr{O}) \leq Ma_{BD}(\mathscr{O})$ with $Ma_{BD}$ given by (\ref{funBD}) and for any $\mathscr{O} \in S=\text{CCG}$ or $\mathscr{O} \in S=\text{CPoly}$, $\ell(s(\mathscr{O})) \leq a(\mathscr{O}) \leq Ma_{2PD}(\mathscr{O})$ with $Ma_{2PD}$ given by (\ref{func2pd}), getting in this way estimates for the assembly complexity of colored Strings, Graphs and Polyominoes. Moreover the upper bounds for $j$-Strings reduce to the well-known bounds derived by Schonhage (\ref{ACineq}) in \cite{SCHONHAGE19751} when $j=1$. \\\\
Lower and upper bounds for the Assembly Index of any Assembly Space $(S,\circ,BB)$ could be used to derive upper and lower bounds for the \textit{complexity level sets} of the Assembly Spaces (i.e. the sets of objects with the same Assembly Index), respectively. The results presented in this paper allow to derive numerical evidence regarding the growth of the combinatorial Assembly Universes as conjectured in \cite{Sharma2023AssemblyTE}, at least for the sets of Strings, Graphs, Polyominoes and Molecules. This is currently work in progress \cite{boundschemistry}. \\\\
Since objects characterized by Optimal Assembly Addition Chains for which every object is obtained by gluing the predecessor object with itself, e.g. the \textit{monochromatic string-like} objects, are intuitively expected to be more symmetrical than objects with the same size characterized by Optimal Assembly Addition Chains where every object is the gluing of two different predecessor objects, it is natural to conjecture that there must be a correlation between the Assembly Index of an object $\mathscr{O}$ and its symmetry \cite{johnston2022symmetry}, \cite{vimal2025open}  (e.g. the length of the longest repeated sub-string, the order of its automorphism group, etc) and checking this is worth of further investigation.

\section*{Appendix A}\label{sec:AppA} 

Let $(S,\circ,BB)$ be an Assembly Space and $\mathscr{O} \in S$. For every $C \in AAC(\mathscr{O},BB)$ we could define a Directed Acyclic Graph $G(C)$ with: 
\begin{enumerate}[(1.)]
    \item At most $s(\mathscr{O})$ sources (i.e. the vertices with indegree equal to zero), each corresponding to the elements $BB$ in $\mathscr{O}$. 
    \item One sink (i.e. only one vertex with outdegree equal to zero), which corresponds $\mathscr{O}$.
    \item $L(C)$ vertices of with indegree equal to $2$ and outdegree of at least one, which correspond to the elements of $C$.
\end{enumerate}
such that there is a directed edge from the vertex corresponding to the object $\mathscr{P} \in C \cup BB$ to the vertex corresponding to the object $\mathscr{Q} \in C$ iff exits $\mathscr{U} \in C \cup BB$ such that
\begin{equation*}
    \mathscr{Q} \in \{\mathscr{P}\} \circ \{\mathscr{U}\}.
\end{equation*}
The map $G: AAC(\mathscr{O},BB) \to DAG$ is well-defined for any $\mathscr{O} \in S$ and if we denote by 
\begin{equation}
    e_+: DAG(\mathcal{S}; BB) \to \mathbb{Z}^+
\end{equation}
the function that assign to each DAG the number of vertices with positive indegree (i.e. ignoring the sources), then the Assembly Index of the object $\mathscr{O} \in S$ will be given by 
\begin{equation}\label{AI}
    a(\mathscr{O}) = \min(e_+(G(AAC(\mathscr{O},BB)))).
\end{equation}
This is the mathematical formalization of the definition of Assembly Index introduced in \cite{Marshall2022FormalisingTP}.

\printbibliography 

@article{johnston2022symmetry,
  title={Symmetry and simplicity spontaneously emerge from the algorithmic nature of evolution},
  author={Johnston, Iain G and Dingle, Kamaludin and Greenbury, Sam F and Camargo, Chico Q and Doye, Jonathan PK and Ahnert, Sebastian E and Louis, Ard A},
  journal={Proceedings of the National Academy of Sciences},
  volume={119},
  number={11},
  pages={e2113883119},
  year={2022},
  publisher={National Academy of Sciences}
}

@article{THURBER2021112200,
title = {Addition chains, vector chains, and efficient computation},
journal = {Discrete Mathematics},
volume = {344},
number = {2},
pages = {112200},
year = {2021},
issn = {0012-365X},
doi = {https://doi.org/10.1016/j.disc.2020.112200},
url = {https://www.sciencedirect.com/science/article/pii/S0012365X20303861},
author = {Edward G. Thurber and Neill M. Clift}
}

@article{brau,
author = {Alfred Brauer},
title = {{On addition chains}},
volume = {45},
journal = {Bulletin of the American Mathematical Society},
number = {10},
publisher = {American Mathematical Society},
pages = {736 -- 739},
year = {1939},
}

@article{erdos1960remarks,
  title={Remarks on number theory III. On addition chains},
  author={Erd{\"o}s, Paul},
  journal={Acta Arithmetica},
  volume={6},
  pages={77--81},
  year={1960},
  publisher={Instytut Matematyczny Polskiej Akademii Nauk}
}

@article{SCHONHAGE19751,
title = {A lower bound for the length of addition chains},
journal = {Theoretical Computer Science},
volume = {1},
number = {1},
pages = {1-12},
year = {1975},
issn = {0304-3975},
doi = {https://doi.org/10.1016/0304-3975(75)90008-0},
url = {https://www.sciencedirect.com/science/article/pii/0304397575900080},
author = {Arnold Schönhage}
}

@article{NPAC,
author = {Downey, Peter and Leong, Benton and Sethi, Ravi},
title = {Computing Sequences with Addition Chains},
journal = {SIAM Journal on Computing},
volume = {10},
number = {3},
pages = {638-646},
year = {1981},
doi = {10.1137/0210047},
URL = {https://doi.org/10.1137/0210047},
eprint = {https://doi.org/10.1137/0210047}
}

@article{kempes2024assembly,
author={Kempes, Christopher P.
and Lachmann, Michael
and Iannaccone, Andrew
and Matthew Fricke, G.
and Redwan Chowdhury, M.
and Walker, Sara I.
and Cronin, Leroy},
title={Assembly theory and its relationship with computational complexity},
journal={npj Complexity},
year={2025},
month={Sep},
day={03},
volume={2},
number={27},
URL={https://doi.org/10.1038/s44260-025-00049-9}
}

@article{boundschemistry,
author={Morales, Juan C.
and Patarroyo, Keith Y.
and Sharma, Abhishek
and Alobo, David O.
and Cronin, Leroy},
title={Bounds of Assembly Theory and the Rate of Growth of Chemical Space},
journal={In preparation},
year={2025},
}

@book{golomb,
  title={Polyominoes: Puzzles, Patterns, Problems, and Packings},
  author={Golomb, S. W.},
  publisher={Princeton University Press},
  year={1995},
edition={2nd ed.}
}

@article{REDELMEIER1981191,
title = {Counting polyominoes: Yet another attack},
journal = {Discrete Mathematics},
volume = {36},
number = {2},
pages = {191-203},
year = {1981},
issn = {0012-365X},
doi = {https://doi.org/10.1016/0012-365X(81)90237-5},
url = {https://www.sciencedirect.com/science/article/pii/0012365X81902375},
author = {D.Hugh Redelmeier}
}

@book{bergman,
author = {Bergman, C.},
year = {2011},
month = {01},
pages = {1-299},
title = {Universal algebra: Fundamentals and selected topics},
journal = {Universal Algebra: Fundamentals and Selected Topics}
}

@article{altman2015integer,
  title={Integer complexity and well-ordering},
  author={Altman, Harry},
  journal={Michigan Mathematical Journal},
  volume={64},
  number={3},
  pages={509--538},
  year={2015},
  publisher={University of Michigan, Department of Mathematics}
}

@article{Stolarsky_1969,
title={A Lower Bound for the Scholz-Brauer Problem}, 
volume={21}, 
DOI={10.4153/CJM-1969-077-x},
journal={Canadian Journal of Mathematics},
author={B. Stolarsky, Kenneth}, 
year={1969}, 
pages={675–683}}

@article{Sharma2023AssemblyTE,
  title={Assembly theory explains and quantifies selection and evolution},
  author={Abhishek Sharma and D{\'a}niel Cz{\'e}gel and Michael Lachmann and Christopher P. Kempes and Sara I. Walker and Leroy Cronin},
  journal={Nature},
  year={2023},
  volume={622},
  pages={321 - 328},
  url={https://api.semanticscholar.org/CorpusID:263669990}
}

@article{Marshall2022FormalisingTP,
  title={Formalising the Pathways to Life Using Assembly Spaces},
  author={Stuart M. Marshall and Douglas G. Moore and Alastair R. G. Murray and Sara Imari Walker and Leroy Cronin},
  journal={Entropy},
  year={2022},
  volume={24},
  url={https://api.semanticscholar.org/CorpusID:250132473}
}

@article{LehmanPhD,
author = {Lehman, Eric},
year = {2002},
month = {03},
pages = {},
journal={Massachusetts Institute of Technology, Dept. of Electrical Engineering and Computer Science, PhD. Thesis},
title = {Approximation Algorithms for Grammar-Based Data Compression}
}

@Article{math12101600,
AUTHOR = {Łukaszyk, Szymon and Bieniawski, Wawrzyniec},
TITLE = {Assembly Theory of Binary Messages},
JOURNAL = {Mathematics},
VOLUME = {12},
YEAR = {2024},
NUMBER = {10},
ARTICLE-NUMBER = {1600},
URL = {https://www.mdpi.com/2227-7390/12/10/1600},
ISSN = {2227-7390},
DOI = {10.3390/math12101600}
}

@article{Scholz,
  title={Aufgabe 253},
  author={Scholz, A},
  journal={Jahresbericht der Deutschen Mathematiker-Vereingung},
  year={1937},
  volume={47},
  pages={41-42},
}

@article{seet2024rapid,
author = {Seet, Ian and Patarroyo, Keith Y. and Siebert, Gage and Walker, Sara I. and Cronin, Leroy},
title = {Rapid Exploration of the Assembly Chemical Space of Molecular Graphs},
journal = {Journal of Chemical Information and Modeling},
year = {2025},
doi = {10.1021/acs.jcim.5c01964},
note ={PMID: 41353721},
URL = {https://doi.org/10.1021/acs.jcim.5c01964},
eprint = {https://doi.org/10.1021/acs.jcim.5c01964}}

@inproceedings{Lehman,
author = {Lehman, Eric and Shelat, Abhi},
title = {Approximation algorithms for grammar-based compression},
year = {2002},
isbn = {089871513X},
publisher = {Society for Industrial and Applied Mathematics},
address = {USA},
booktitle = {Proceedings of the Thirteenth Annual ACM-SIAM Symposium on Discrete Algorithms},
pages = {205–212},
numpages = {8},
location = {San Francisco, California},
series = {SODA '02}
}

@article{bondy1990perfect,
  title={Perfect path double covers of graphs},
  author={Bondy, J Adrian},
  journal={Journal of graph theory},
  volume={14},
  number={2},
  pages={259--272},
  year={1990},
  publisher={Wiley Online Library}
}

@article{clift,
author = {Clift, Neill Michael},
title = {Calculating optimal addition chains},
year = {2011},
issue_date = {March 2011},
publisher = {Springer-Verlag},
address = {Berlin, Heidelberg},
volume = {91},
number = {3},
issn = {0010-485X},
url = {https://doi.org/10.1007/s00607-010-0118-8},
doi = {10.1007/s00607-010-0118-8},
journal = {Computing},
month = mar,
pages = {265–284},
numpages = {20}
}

@article{kotzig1957theorie,
  title={Aus der Theorie der endlichen regul{\"a}ren Graphen dritten und vierten Grades},
  author={Kotzig, Anton},
  journal={Casopis Pest. Mat},
  volume={82},
  pages={76--92},
  year={1957}
}

@article{RSPaper,
    author = {Marshall, Stuart M. and Murray, Alastair R. G. and Cronin, Leroy},
    title = {A probabilistic framework for identifying biosignatures using Pathway Complexity},
    journal = {Philosophical Transactions of the Royal Society A: Mathematical, Physical and Engineering Sciences},
    volume = {375},
    number = {2109},
    pages = {20160342},
    year = {2017},
    month = {11},
    issn = {1364-503X},
    doi = {10.1098/rsta.2016.0342},
    url = {https://doi.org/10.1098/rsta.2016.0342},
    eprint = {https://royalsocietypublishing.org/rsta/article-pdf/doi/10.1098/rsta.2016.0342/389700/rsta.2016.0342.pdf},
}

@article{vimal2025open,
  title={Open, Reproducible Calculation of Assembly Indices},
  author={Vimal, Devansh and Parzych, Garrett and Smith, Olivia M and Parkar, Devendra and Bergen, Sean and Daymude, Joshua J and Mathis, Cole},
  journal={arXiv preprint arXiv:2507.08852},
  year={2025}
}

\end{document}